\documentclass [11pt]{article}
\usepackage[utf8]{inputenc} 
\usepackage[T1]{fontenc} 
\usepackage{amsmath}
\usepackage{amsfonts}
\usepackage{amsthm}

\usepackage[dvipsnames]{xcolor}

\newcommand{\footnoteremember}[2]{
\footnote{#2}
\newcounter{#1}
\setcounter{#1}{\value{footnote}}
}

\newlength{\oldparindent}
\newcommand{\cL}{{\mathbb {L}}}%espace L^p

\newcommand{\bpf}{\begin{preuve}}
\newcommand{\epf}{ \end{preuve} \medskip}

\newcommand{\benum}{\begin{enumerate}}
\newcommand{\eenum}{\end{enumerate}}

\newcommand{\bitem}{\begin{itemize}}
\newcommand{\eitem}{\end{itemize}}

\newcommand{\brmq}{\begin{rmq}}
\newcommand{\ermq}{\end{rmq}}

\newcommand{\brmqs}{\begin{rmqs}}
\newcommand{\ermqs}{\end{rmqs}}

\newcommand{\bapp}{\begin{application}}
\newcommand{\eapp}{\end{application}}

\newcommand{\bapps}{\begin{applications}}
\newcommand{\eapps}{\end{applications}}

\newcommand{\bdefi}{\begin{definition}}
\newcommand{\edefi}{\end{definition}}

\newcommand{\beq}{\begin{equation}}
\newcommand{\eeq}{\end{equation}}

\def\bpm{\begin{pmatrix}}
\def\epm{\end{pmatrix}}

\newcommand{\bcas}{\begin{cases}}
\newcommand{\ecas}{\end{cases}}

\newcommand{\bex}{\begin{exemp}}
\newcommand{\eex}{\end{exemp}}

\newcommand{\bexs}{\begin{exemps}}
\newcommand{\eexs}{\end{exemps}}

\newcommand{\bprop}{\begin{proposition}}
\newcommand{\eprop}{\end{proposition}}

\newcommand{\bthm}{\begin{theoreme}}
\newcommand{\ethm}{\end{theoreme}}

\newcommand{\bcor}{\begin{corollaire}}
\newcommand{\ecor}{\end{corollaire}}

\newcommand{\blem}{\begin{lemme}}
\newcommand{\elem}{\end{lemme}}

\newcommand{\beqna}{\begin{eqnarray}}
\newcommand{\eeqna}{\end{eqnarray}}

\newcommand{\beqnas}{\begin{eqnarray*}}
\newcommand{\eeqnas}{\end{eqnarray*}}

\newcommand{\cA}{{\mathcal A}}

\definecolor{green}{rgb}{0,.7,.2}
\definecolor{orange}{rgb}{0.9,.5,0}

\newcommand{\LL}{{\rm L}}%generator
\def\I{1\!{\rm l}}
\def\E{{\mathbb{E}}}%esperance
\def\tr{{\rm trace}\,}

\def\lag{\langle}
\def\rag{\rangle}

\def\det{{ \rm{det}}}  %parentheses have been corrected
\def\Id{{\rm{Id}}} %parentheses have been corrected

\def\cA{{\mathcal A }}

\def\cD{{\mathcal D}}

\def\cH{{\mathcal  H}}

\def\cL{{\mathcal L }}
\def\cM{{\mathcal M }}
\def\cP{{\mathcal P }}

\newcommand{\bB}{{\mathbb {B}}}
\newcommand{\bC}{{\mathbb {C}}}%complexes

\newcommand{\bE}{{\mathbb {E}}}

\newcommand{\bH}{{\mathbb {H}}}%quaternions , espaces hyperbolique

\newcommand{\bN}{{\mathbb {N}}}% entiers

\newcommand{\bR}{{\mathbb {R}}}%reels
\newcommand{\bS}{{\mathbb {S}}}

\newtheorem{theoreme}{Theorem}[section]
\newtheorem{lemme}[theoreme]{Lemma}
\newtheorem{definition}[theoreme]{Definition}
\newtheorem{proposition}[theoreme]{Proposition}
\newtheorem{corollaire}[theoreme]{Corollary}

\newtheorem{defi}[theoreme]{Definition}
\newenvironment{exemp}{\noindent{\bf Example. --- }}{\par}
\newenvironment{exemps}{\noindent{\bf Examples}\benum}{\eenum\par}
\newtheorem{rmq}[theoreme]{Remark}
\newtheorem{rmqs}[theoreme]{Remarks}

\newenvironment{preuve}{\noindent{\it Proof. --- }}
{\hfill\rule{1.3mm}{2mm}\par} 
\newenvironment{application}{\noindent{\bf Application. --- }}{\par}
\newenvironment{applications}{\noindent{\bf Applications. --- 
}\benum}{\eenum\par}

\newcommand{\discr}{{\rm discr} }
\newcommand{\DS}{\mathcal {\Delta}}
\newcommand{\Cl}{\mathcal {C}l}

\begin{document}

\author{D. Bakry \footnoteremember{IMT}{Institut de Math\'ematiques de Toulouse, Universit\'e Paul Sabatier, 118 
route de Narbonne, 31062 Toulouse, France},  M. Zani \footnote{Laboratoire d'Analyse et Math\'ematiques appliqu\'ees , UMR CNRS 8050, Universit\'e Paris-Est -Cr\'eteil, 61 av du G\'en\'eral de Gaulle, 94010, Cr\'eteil Cedex, France}}
\date{\today}

\title{Random symmetric matrices on Clifford algebras}
\maketitle

\abstract{We consider Brownian motions and other processes (Ornstein-Uhlenbeck processes, spherical Brownian motions)  on various sets of symmetric matrices constructed from algebra structures, and look at their associated spectral measure processes. This leads to the identification of the multiplicity of the eigenvalues, together with the identification of the spectral measures. For Clifford algebras, we thus recover  Bott's periodicity. }
\section{Introduction}
Many works on random matrix theory deal with the law of the spectrum on some specific sets: real symmetric, Hermitian, orthogonal or unitary, etc.
There exists a large literature on this topic. We shall mention the early works of Wigner \cite{Wig58}, Dyson \cite{dyson}, Mehta \cite{Meh}, Mar\v{c}enko and Pastur \cite{marcpastur}, Soshnikov \cite{Sos99} and more recently Anderson, Guionnet and Zeitouni \cite{AndGuioZeit}, Erd\"os and co \cite{erdosandco1,erdosandco2,erdosandco3}, Forrester \cite{Forrester} and references therein.
One may also consider stochastic processes on these sets of matrices (Euclidean  Brownian motions, Ornstein-Uhlenbeck operators, spherical Brownian motions, Brownian motions on groups, see e.g. \cite{YanDoumerc}): they are diffusion processes which  are reversible under the laws we consider.  One  then consider the stochastic process which is the empirical measure of the spectrum of the matrix.    It turns out that in many situation, these empirical measures are again stochastic diffusion processes, and their reversible measures are the image of the measure on the set of matrices. It may then be an easy  way of computing  those spectral measures, but the study of those spectral processes is by itself a topic of interest.

For, say Gaussian,  symmetric,  Hermitian or quaternionic matrices, when one considers the law of  their eigenvalues $(\lambda_1, \cdots, \lambda_n)$,  ordered  for example as $\lambda_1< \cdots < \lambda_n$, it has a density with respect to the Lebesgue measure  $d\lambda_1\cdots d\lambda_n$ which is $C(\prod_{i<j} |\lambda_i-\lambda_j|)^a e^{-1/2\sum_i \lambda_i^2}$, where $C$ is a normalizing constant and  $a=1,2,4$ according to the fact that we are in the real, complex or quaternionic case.  This factor $(\prod_{i<j} |\lambda_i-\lambda_j|)^a$ also appears in many other situations ($SO(n), SU(n), Sp(n)$ matrices, for example). There are even some results for octonionic $2\times2$ matrices where $a=8$.  On the other hand, if one considers the real symmetric, Hermitian and quaternionic $n\times n$ matrices as real matrices (with size $n\times n, (2n)\times (2n)$ and $(4n)\times (4n)$ respectively), their eigenspaces have dimension $1,2$ and $4$ respectively (and $8$ in the special case of octonionic matrices). Therefore, one may  think that this factor $a$ is due to the multiplicity of the eigenvalues for the real form of the matrices. We shall see  that this is not the case. 

  In this paper,  we go beyond the cases of real symmetric, Hermitian and quaternionic cases, and we consider some associative dimension $2^p$ algebras of the Clifford type. With the help of the algebra structure, one may define real symmetric matrices  with size $(n\times 2^p)\times (n\times 2^p)$, which have eigenspaces with dimension $2^k$, where $k$ is related to $p$ in some specific way described below, and  may be chosen as large as we wish. We then consider Brownian motion on these set of matrices, compute the spectral measure processes, which appear to be in many situations symmetric diffusion processes. This fact depends indeed on the algebra structure,  and this leads to the computation of the law of the spectrum. It turns out that there is still a factor $(\prod_{i<j} |\lambda_i-\lambda_j|)^a$ with $a=1,2,4$, and we do not produce in this way other values for $a$ : these values reflect in fact a phenomenon which is known as Bott's periodicity, and has nothing to to with the dimension of the eigenspaces  ( for references on Bott's Periodicity Theorem, see \cite{atiyahbottshapiro,atiyahsinger,karoubi}
).
  
  In order to deal with spectral measures, we consider processes on the characteristic polynomials $P(X)= \det(M-X\Id)$. Indeed, functions of the spectral measures are nothing else than symmetric functions of the  roots of $P$. Usually, in order to characterize the laws of the spectral measure, one works with  its  moments, that is the functions $\sum_{i}^n \lambda_i^k$, where $(\lambda_1, \cdots, \lambda_n)$ denote the eigenvalues of the matrix. We find more convenient to deal with the elementary symmetric functions of the roots, that is the coefficients of $P$. Curiously, this approach is not that popular, although very close to the study of Stieltjes transform of the measure, and we therefore develop some light machinery  in order to perform those computations. To understand the various quantities appearing in the computations, we first consider the diagonal cases (that is when we start directly on some processes on the roots of $P(X)$), that we analyze in the case of flat Brownian motion, Ornstein-Uhlenbeck operators and spherical brownian motions. Then, we pass to the analysis of the laws of the characteristic polynomials in various sets of symmetric matrices, starting with the classical settings(real symmetric, Hermitian and quaternionic), before going to the case of general Clifford algebras.
  
  The paper is organized as follows.
  Section~\ref{sec.symm.diff} is  a short introduction to the methods and language of symmetric diffusion processes. In  Sections~\ref{sec.eucl.lapl},~\ref{sec.diag.OU} and~\ref{sec.diag.sph}   we consider the simpler  cases of diagonal matrices in the flat Euclidean, Gaussian and spherical cases respectively.  We  describe how to handle and identify  various quantities  (discriminants, metric structures, etc) which appear in the computations  through the characteristic polynomials.   Sections~\ref{sec.sym.matr} and~\ref{sec.herm.matr} are devoted to the analysis of the real symmetric, Hermitian and quaternionic cases (although the quaternionic case is just sketched since it is not simpler to handle than the general Clifford case). In Section~\ref{sec.gal.cliff} we introduce the general Clifford algebras we are dealing with, through a presentation which is quite handy  for our purpose and does not seem to be classical. Finally, in Section~\ref{sec.stand.cliff}, we give the complete description of the laws of the spectra for the standard Clifford algebras, where we recover Bott's periodicity through the various laws appearing in the spectra of the matrices , for different values of  the dimension of the algebra. For references on Clifford algebras see \cite{lounesto} and for  some generalized Clifford algebras see  \cite{RamJag}.

\section{ Symmetric diffusion generators and their images\label{sec.symm.diff}}

The general setting for symmetric diffusion generators we deal with is inspired from \cite{bglbook} and  is the following: let $E$ be a smooth manifold endowed with a $\sigma$-finite measure $\mu$  and let $\cA_0$ be be some algebra of smooth function on $E$;  typically when $E$ is an open set in $\bR^n$, $\cA_0$ is the smooth compactly supported functions or polynomials functions on $E$. 
Any function of $\cA_0$  belongs  to every $\cL^p(\mu)$, for any $1\leq p<\infty$, and  $\cA_0$ is dense in very $\cL^p(\mu)$. Moreover,  $\cA_0$ is stable under any transformation $(f_1, \cdots, f_n)\mapsto \Phi(f_1, \cdots, f_n)$, where $\Phi$ is a smooth function $\bR^n\mapsto \bR$ such that $\Phi(0)=0$.  For any linear operator $\LL : \cA_0\mapsto \cA_0$, one defines its carré du champ operator
$$\Gamma(f,g) = \frac{1}{2}\Big( \LL(fg)-f\LL(g)-g\LL(f)\Big)\,.$$
We have the following
\begin{defi}\label{defi.diff}
A symmetric diffusion operator is a  linear operator $\LL$: $\cA_0\oplus 1\mapsto \cA_0$, such that 
\benum
\item $ \LL(1)=0$,
\item $\forall f,g\in \cA_0\oplus 1, ~\int f\LL(g)\, d\mu = \int g \LL (f) \, d\mu$,
\item $\forall f\in \cA_0, \Gamma(f,f) \geq 0$,
\item\label{change.variable} $\forall f=(f_1, \cdots, f_n)$, where $f_i\in \cA_0$ , $\forall \Phi : ~ \bR^n\mapsto \bR$, 
$$\LL(\Phi(f))= \sum_i \partial_i \Phi(f) \LL(f_i) + \sum_{i,j} \partial^2_{ij} \Phi(f) \Gamma(f_i,f_j).$$
\eenum
\end{defi}
Such an operator describes the law of a stochastic process with generator $\LL$ and reversible measure $\mu$. It is often important to be able to extend $\LL$ on a larger class of functions (typically smooth functions), and we shall do this without further comments.  For a consistent reference on these operators, see \cite{bglbook}.

Let us consider an open set $\Omega\subset E$, and a given system of coordinates $(x^i)$, then we can write
$$\LL(f) = \sum_{ij} g^{ij} (x)\partial^2_{ij} f + b^i(x)\sum_i \partial_i f,$$ where 
$$g^{ij}(c)= \Gamma(x^i, x^j), ~b^i (x)= \LL(x^i).$$
 
The positivity condition (see point 3. in Definition~\ref{defi.diff}  above)  imposes that for any $x\in \Omega$, the symmetric matrix $(g^{ij}(x))$ is non negative. Moreover, provided $\mu$ has a smooth positive density $\rho$ with respect to  the Lebesgue measure $dx^1\cdots dx^n$, the operator $L$ is entirely described, up to a normalizing factor, by 
\beq\label{eq.rho}\LL(f)= \frac{1}{\rho}\sum_{ij} \partial_j(g^{ij}\rho\partial_j f).\eeq
Suppose that one may find functions $(a^1, \cdots, a^k)$ such that for some smooth functions $B^i$ and $G^{ij}$
$$\LL(a^i)= B^i(a^1, \cdots, a^k), ~\Gamma(a^i, a^j)= G^{ij}(a^1, \cdots, a^k),$$ then, writing $a= (a^1, \cdots, a^k)$, we get readily
$$\LL (f(a))= \sum_{ij} G^{ij}(a) \partial^2_{ij} f (a) + \sum_i B^i(a)\partial_i f (a),$$ and  the operator 
$$\hat \LL= \sum_{ij} G^{ij} \partial^2_{ij} + \sum_iB^i\partial_i$$ is nothing else than the operator $\LL$ acting on functions depending only on $a=(a^1, \cdots, a^k)$. We shall say that $\hat \LL$ is the image of the operator $\LL$ through $a=(a^1, \cdots, a^k)$. As a consequence, the operator $\hat \LL$ is  symmetric with respect to the image measure of $\mu$ through $a=(a^1, \cdots, a^k)$. With the help of equation~\eqref{eq.rho}, it will be a good way to identify the image measure.

A special case concerns Laplace operators, where $(g^{ij})$ is everywhere non degenerate and $\rho$ is by definition $\det( g)^{-1/2}$. When $x\mapsto a=(a^1, \cdots, a^n)$ is a local diffeomorphism (a simple change of coordinates), then the image of $\LL$ is again a Laplace operator, (the Laplace operator written in the new system of coordinates), and the density measure is again $\det(G)^{-1/2}$ is the new system of coordinates. It is however important to notice that we shall use this procedure even when $x\mapsto a$ is not a diffeomorphism, for example with the map $M\mapsto P(X)$, where $M$ is a matrix and $P$ its characteristic polynomial.

Those operators $L$ are related with the associated stochastic processes  $(X_t)$ by the requirement that, for any smooth function $f$, $f(X_t)-\int_0^t \LL(f)(X_s) \,ds$ is a local martingale (a true martingale if for example $f$ and $L(f)$ are bounded).  This is enough to describe  the law of $(X_t)$, or the joint laws of $(X_{t_1}, \cdots, X_{t_n})$ from the starting point $X_0=x$ and the knowledge of the operator $\LL$, see~\cite{bglbook}. Anyhow, we shall not really use this interpretation in terms of stochastic processes in what follows, since we shall mainly concentrate on the properties of the operator  $\LL$.
As mentioned above, we deal with operators acting on the space of polynomials i.e. on $\bR^n$, when identifying $\bR^n$ with the set of  the coefficients of the polynomial. The polynomials shall be monic in general, i.e.  $P(X)= X^n + \sum_0^{n-1} a_i X^i$, where $(a_0, \cdots, a_{n-1})$ is the stochastic process, and $X$ may be considered as a parameter. 
The coefficients $(a_i)$ can be viewed as coordinates in this set of polynomials, writing for example $\LL\big(P(X)\big)= \sum_i X^i \LL(a_i)$, and doing similarly for the operator $\Gamma\big(P(X), P(Y)\big)$. One may also consider a fixed $X$ and see $P(X)$ as an application $\bR^n \mapsto \bR$, for which one can take $\log (P(X))$ or $P(X)^\beta$. 
Last, those function can be described as series in the variable $X$ ( even formal series, regardless of their domain of convergence), all computations on these expressions boiling down to algebraic computations involving $\LL$,  $\Gamma$ and polynomial expressions in the coefficients $a_i$ (see for example Lemma~\ref{lem.puiss.P} below).

\section{ The image of the Euclidean Laplacian under elementary symmetric functions\label{sec.eucl.lapl}}
Let $x= (x_1, \cdots , x_n)\in \mathbb R^n$ and 
$$P(X) = \prod_{i=1}^n  (X-x_i)= \sum_{i= 0}^n a_iX^i,$$ such that $ (-1)^ia_i(x_1, \cdots, x_n)$ are the elementary symmetric functions.
If we want to describe the image of the Laplace operator $\Delta_{\bE}$ on $\mathbb R^n$ under
symmetric functions of $(x_1, \cdots, x_n)$, we  may look at  smooth functions $F(a_0, \cdots, a_{n-1})$. At least in the Weyl chamber $\{x_1< x_2 \cdots < x_n\}$ the application 
$(x_1, \cdots, x_n)\mapsto \Phi(x_1, \cdots, x_n)= (a_0, \cdots, a_{n-1})$ is a local diffeomorphism.
We first have to look at the image of the Lebesgue measure $dx= dx_1\cdots dx_n$ under $\Phi$. For that purpose, we recall the definition of the discriminant of a polynomial $P$.
We consider monic polynomials: for two monic  polynomials $P(X)= \sum_{i=0}^n a_i X^i$ and $Q(X)= \sum_{i=0}^p b_i X^i$, the resultant 
$R(P,Q)$ is a polynomial in the coefficients $(a_0, \cdots, a_{n-1}, b_0, \cdots, b_{p-1})$ which vanishes exactly when $P$ and $Q$ have a common root (in the complex plane).
Indeed, $R(P,Q)$ is the determinant of the $n\times p$ Sylvester matrix
$$\bpm 	1& a_{n-1} & a_{n-2}&\cdots & a_0&0&\cdots&0\\
		0&1         	& a_{n-1}	&\cdots & a_1&a_0&\cdots &0\\
		0&0		&1        	 & \cdots     &a_2&a_1&\cdots & 0\\
		\cdots & \cdots &\cdots                   &\cdots           & \cdots      &  \cdots     & \cdots&\cdots\\
		\cdots & \cdots &     \cdots              &      \cdots     &    a_{p-2}   &  \cdots     & a_1&a_0\\
		1& b_{p-1} & b_{p-2}&\cdots & b_0&0&\cdots&0\\
		0&1         	& b_{p-1}	&\cdots & b_1&b_0&\cdots &0\\
		0&0		&1        	 & \cdots     &b_2&b_1&\cdots & 0\\
		\cdots & \cdots &\cdots                   &\cdots           & \cdots      &  \cdots     & \cdots&\cdots\\
		\cdots & \cdots  &      \cdots              &      \cdots      &     \cdots   &  \cdots     & b_1&b_0

\epm$$
It can be viewed as the determinant of the following system of linear equations in the unknown variables $\{1,X,\cdots, X^{n+p-1}\}$:
$$\{P(X)=0, XP(X)=0, \cdots, X^{p-1}P(X)=0, Q(X)=0, XQ(X)=0, \cdots, X^{n-1} Q(X)=0\}$$ 
It turns out that, when $P(X)= \prod_i(X-x_i)$ and $Q(X)= \prod_j(X-y_j)$, then 
$R(P,Q)= \prod_{i,j} (x_i-y_i)$.

The discriminant $\discr(P)$ is $D(P)= (-1)^{n(n-1)/2}R(P,P')$ and expresses a necessary and sufficient condition for $P$ to have a double root.
Then, when $P(X)= \prod(X-x_i)$, one has 
$\discr(P)= \prod_{i<j} (x_i-x_j)^2$.

The discriminant is not an easy expression of the coefficients $(a_0, \cdots, a_{n-1})$, and the following computations are here to make them easier.

\bprop The image measure of $dx$ under $\Phi$ is $$d\mu_0=n! |\discr(P)|^{-1/2}\I_{\cD>0} da_0\cdots da_{n-1},$$where $D$ is the connected component
of the set $\{\discr(P)>0\}$ where all the roots of the polynomial $P$ are real.

\eprop
\bpf 

Let is first observe that for a polynomial $P$ having only distinct real roots, $\discr(P)>0$. However, the condition $\discr(P)>0$  is not sufficient to assert that $P$ has only real roots: it is also positive  for example when $P$ has an even number of pairs of complex conjugate roots. The set where all the roots are real and distinct is obviously connected, and therefore there is only one connected component of this set where all the roots are real. This is our set $\cD$, the set which contains for example the point $\prod_0^{n-1}(X-i)$.

To compute the image measure of $\Phi$, it is enough to identify the image of $\Phi$ when restricted to the Weyl chamber $\{x_1<\cdots < x_n\}$, since the same computation will hold true in any other Weyl chamber (that is the image of the first one under a permutation of the coordinates), and there are $n!$ of such chambers.

 First, it is clear that the support of the image measure is included in the closure of  $\cD\subset \{\discr(P)>0\}$, and, from the explicit expression of the discriminant in terms of the roots, is  strictly positive  in any Weyl chamber.  The boundary of this set is a subset of the algebraic surface (in the space of the $(a_i)$ coordinates) $\{ \discr(P)=0\}$.

Indeed, it is quite easy to see this result by induction on the degree $n$. It is clear that it is true for $n=1$, since $a_1= -x_1$

Let us assume that the result is true for $n$ and set 
$P(X)= (X-x_{n+1})Q(x)= \sum_0^{n+1} a_i X^i$, where 
$Q(X)= \sum_{0}^n b_i X^i= \prod_1^{n}(X-x_i)$. Then
$$a_{n}= b_{n-1}-x_{n+1}, a_{n-1}= b_{n-2}-x_{n+1}b_{n-1}, \cdots, a_1= b_{0}-x_{n+1} b_{1}, a_{0}= -x_{n+1}b_{0}.$$
The Jacobian of the transformation 
$(x_{n+1}, b_0, b_1, \cdots, b_{n-1})\mapsto (a_{0}, \cdots, a_n)$ is easily seen to be 
$$|x_{n+1}^{n}+ b_nx_{n+1}^{n-1} + \cdots + b_0|= |Q(x_{n+1})|= |\prod_1^{n} (x_{n+1}-x_i)|.$$
Therefore, if 
$db_0\cdots db_{n-1}= |\prod_{1\leq i<j\leq n} (x_i-x_j)| dx_1\cdots dx_{n}$, then 
$$da_0\cdots da_n= |\prod_{1\leq i<j\leq n+1}(x_i-x_j)|.$$

From what precedes, it is clear that the Jacobian of the transformation $(x_1, \cdots, x_n)\mapsto (a_0, \cdots, a_{n-1})$ is non zero on any Weyl chamber, and then the boundary of the image is included in the algebraic set $ \{\discr(P)=0\}$. This is enough to identify the support of the image measure as the closure of $\cD$.

\epf

Let us now compute the image of the Laplace operator $\Delta_{\bE}$ on $\bR^n$ under $\Phi$. In what follows,  and throughout the paper, $\Gamma_\bE$ denotes the Euclidean carré du champ, that is, in the standard system of coordinates, 
$$\Gamma_\bE(f,g) = \sum_i \partial_i f\partial_i g.$$

Fix $X \in \bR$, and consider the function
\begin{eqnarray*}
\bR^n\setminus \{ \exists i, x_i=X\}&\mapsto& \bR \\
(x_1, \cdots, x_n) &\mapsto& \log P(X)= \sum_i \log (X-x_i)
\end{eqnarray*} We get

\bprop
\benum
\item For any $X\in \bR$, 
$\Delta_{\bE}(P(X))=0$
\item
For any $(X,Y)\in \bR^2$,
\beq\label{eq.gamma.logP}\Gamma_\bE(\log(P(X), \log P(Y)) = \frac{1}{Y-X}\Big(\frac{P'(X)}{P(X)}-\frac{P'(Y)}{P(Y)}\Big).\eeq
\eenum

\eprop

\bpf The first assertion is immediate, since every function $a_i$ is an harmonic function on $\bR^n$ ( as a polynomial of degree 1 in any coordinate $x_i$).

For the second one, one has
$$\Gamma_\bE(\log P(X), \log P(Y))= \sum_{i} \partial_{x_i} \log P(X) \partial_{x_i} \log P(Y)= \sum_i \frac{1}{(X-x_i)(Y-x_i)}.$$
But 
$$\frac{1}{(X-x_i)(Y-x_i)}= \frac{1}{Y-X}(\frac{1}{X-x_i}-\frac{1}{Y-x_i})$$ and 
$$\sum_i\frac{1}{(X-x_i)(Y-x_i)}= \frac{1}{Y-X} \Big(\frac{P'(X)}{P(X)}-\frac{P'(Y)}{P(Y)}\Big).$$
\epf

\bcor\label{cor.val.Gamma} Setting  $\alpha_{i,j}= (i+1)a_{i+1}a_j$, where  $a_i= 0 $ if $i>n$ and  $a_n=1$, one has 
$$\sum_{i,j} \Gamma_\bE(a_i,a_j)X^iY^j = \sum_{i\neq j }\alpha_{i,j} \frac{X^iY^j-X^jY^i}{Y-X},$$  from which
\beq\label{eq.deriv.G}\Gamma_\bE(a_k,a_p)= \sum_{(p-k)_+\leq l\leq p} \alpha_{p-l, k+l+1}-\sum_{(k-p)_+\leq l \leq k} \alpha_{p+l+1,k-l}.\eeq

Moreover
\beq\label{eq.deriv.G1}\partial_{a_k}\Gamma_\bE(a_k, a_p)= \I_{p\geq k} (k-p-2) a_{p+2}+ \I_{p= k-1} ka_{p+2}.\eeq
and 
\beq\label{eq.deriv.G2}\sum_k \partial_{a_k} \Gamma_\bE(a_k, a_p) = -\frac{1}{2}(p+1)(p+2) a_{p+2}.\eeq

\ecor

\bpf This is straightforward from equation~\eqref{eq.gamma.logP}, which gives
$$\Gamma_\bE(P(X), P(Y))= \frac{1}{Y-X}(P'(X)P(Y)-P'(Y)P(X)).$$

On the other hand, by bilinearity, one has
$$\Gamma_\bE(P(X),P(Y))= \sum_{ij} \Gamma_\bE(a_i,a_j) X^iY^j.$$ 
We then obtain
$$\sum_{i,j} \Gamma_\bE(a_i,a_j)X^iY^j = \sum_{i\neq j }\alpha_{i,j} \frac{X^iY^j-X^jY^i}{Y-X},$$ 
from which  formula~\eqref{eq.deriv.G} follows easily. 

Formulae~\eqref{eq.deriv.G1} and ~\eqref{eq.deriv.G2} are easy consequences of the first one.

\epf

The application $(x_1, \cdots, x_n)\mapsto (a_0, \cdots, a_{n-1})$ is a local diffeomorphism. Therefore, the image of the Laplace operator in coordinates $(x_1, \cdots, x_n)$  is the Laplace operator in coordinates $(a_0, \cdots, a_{n-1})$. This leads to some formulae which are not immediate.

\bcor\label{cor.gamma.discr} One has
\benum
\item\label{cor.gamma.discr1}  $\discr(P)= \det(\Gamma_\bE(a_i,a_j))$.

\item \label{cor.gamma.discr2} For any $i\in\{0, \cdots, n-1\}$, $\sum_j\Gamma_\bE(a_i,a_j) \partial_{a_j}\log \discr(P)= 2\sum_j \partial_{a_j}\Gamma_{\bE}(a_i,a_j)$.

\item\label{cor.gamma.discr3} For any $i\in\{0, \cdots, n-1\}$, $\sum_{i,j} X^i\Gamma_\bE(a_i,a_j) \partial_{a_j}\log \discr(P)= -P''(X)$.

\eenum

\ecor

\bpf

 One the one hand, we know that in the coordinates $(a_0, \cdots, a_{n-1})$, the Laplace operator has carré du champ $\Gamma_\bE(a_i,a_j)$ and  reversible measure (here the Riemann measure) density $\discr(P)^{-1/2}da_0\cdots da_{n-1}$.

One the other hand, we know that this density measure is always, (for any Laplace operator) $\det (\Gamma_\bE(a_i,a_j)^{-1/2})$. This is enough to get point ~\ref{cor.gamma.discr1}.

Point~\ref{cor.gamma.discr2} comes from the observation that in these coordinates $(a_i)$, $\Delta(a_i)=0$.
Setting $G^{ij}= \Gamma_\bE(a_i,a_j)$, and $\rho= \det(G^{ij})^{-1/2}$, the Laplace operator writes
$$\Delta(f) = \sum_{ij} G^{ij} \partial^2_{ij} f+ \sum_{ij} \partial_{j} (G^{ij}) \partial_i f +\sum_{ij}G^{ij} \partial_i f\partial_j \log \rho.$$
From this, we know that 
$$\Delta(a_i)= \sum_{j} \partial_{j} G^{ij}  -\frac{1}{2}\sum_j G^{ij} \partial_j \log \discr(P).$$ Applying $\Delta(a_i)=0$, point~\ref{cor.gamma.discr2} is the direct translation of the previous.

For point~\ref{cor.gamma.discr3}, it suffices to combine point~\ref{cor.gamma.discr2} together with formula~\eqref{eq.deriv.G} to obtain, for $i= 0, \cdots, n-1$

$$\sum_j\Gamma_\bE(a_i,a_j) \partial_{a_j}\log \discr(P)= -(i+2)(i+1)a_{i+2}.$$

\epf
As a consequence

\bprop\label{cor.gamma.discr4} 

$\Gamma_\bE(P, \log \discr(P))= -P''$.

\eprop

\bpf
This is just a rephrasing of point~\ref{cor.gamma.discr3}. in Corollary~\ref{cor.gamma.discr} .

\epf

\section{The image of the Ornstein-Uhlenbeck operator under elementary symmetric functions\label{sec.diag.OU}}

In Euclidean spaces, the Ornstein-Uhlenbeck operator is defined as
$$\LL_{\mathbb{OU}}(f)= \Delta_\bE (f)- \frac{1}{2}\Gamma_\bE(\|x\|^2, f)= \Delta_\bE(f) -\sum_{i=1}^n x_i \partial_i f.$$ 
It shares the same carré du champ operator  than the Laplace operator, but has the standard Gaussian measure as reversible measure.
It also admits a complete orthonormal system of eigenvectors, namely the Hermite polynomials, which are total degree $k$ polynomials in the variables $(x_1, \cdots, x_n)$ with associated eigenvalue $\LL_{\mathbb{OU}}(H)= -kH$.

\bprop \label{driftOU}For $P(X)= \prod_{i=1}^n(X-x_i)= \sum_{i=0}^n a_iX^i$, one has
\benum
\item\label{driftOU1} $\LL_{\mathbb{OU}}(P)=-\sum_i x_i \partial_i P= \sum_i(n-i)a_i X^i= -nP(X)+XP'(X)$.
\item \label{driftOU2}$\frac{1}{2}\Gamma_\bE(a_{n-1}^2-2a_{n-2}, P(X))= nP(X)-XP'(X)$. 
\item\label{driftOU3}$\forall i= 0, \cdots, n-1$, $a_{n-1}\Gamma_\bE(a_i,a_{n-1})-\Gamma_\bE(a_i, a_{n-2})= (n-i) a_i$.
\item \label{mesureOU} The image of the Gaussian measure  is $$\frac{n!}{(2\pi)^{n/2}} \exp(a_{n-2}-\frac{a_{n-1}^2}{2})\discr(P)^{-1/2}\I_{D} \prod_0^{n-1} da_i= \frac{n!}{(2\pi)^{n/2}} \exp(a_{n-2}-\frac{a_{n-1}^2}{2})d\mu_0.$$
\eenum

\eprop

\bpf  
The first item~\ref{driftOU1} comes from the fact that since $\Delta (a_i)=0$, one has 
$$\LL_{\mathbb{OU}}(P)= \sum_i X^i \LL_{\mathbb{OU}}(a_i)= -\sum_{ij} X^i x_j\partial_j (a_j).$$
But the functions $a_j$ are homogeneous polynomial of degree $n-i$ in the variables $x_i$, and therefore
$\sum_j x_j\partial_j (a_j)= (n-i)a_i$.

If we observe that $\sum_{ij} x_j \partial_i f= \frac{1}{2} \Gamma(\|x\|^2,f)$, and that $\|x\|^2= a_{n-1}^2-2a_{n-2}$, 
one sees that 
$$\frac{1}{2}\Gamma_\bE(a_i, a_{n-1}^2-2a_{n-2})= (n-i) a_i,$$ from which 
$$\frac{1}{2}\Gamma_\bE(a_{n-1}^2-2a_{n-2}, P(X))= nP(X)-XP'(X).$$

Points~\ref{driftOU3}. and~\ref{mesureOU}.  are  immediate consequences.

\epf 

\brmq Observe also that if $U= a_{n-1}^2-2a_{n-2}$, then, for the Euclidean quantities,  $\Gamma_\bE(U,U)= 4U$ and $\Delta_\bE(U)= 2n$, as expected.

\ermq 
It is worth to observe that the measure $\frac{n!}{(2\pi)^{n/2}} \exp(a_{n-2}-\frac{a_{n-1}^2}{2})d\mu_0$, contrary to the Gaussian  measure, does not have in general exponential moments.  Nevertheless, polynomials are dense in $\cL^2(\mu)$. Indeed, any function $f$ in $\cL^2(\mu)$ which is orthogonal  for the Gausian measure to any polynomial $Q(a_0, \cdots, a_{n-1})$ may be lifted into a function $\hat f : \bR^n\mapsto \bR$ which is invariant under permutation of the variables $(x_1, \cdots, x_n)$ and orthogonal to any symmetric polynomial.  But such a function would then be orthogonal  to any polynomial, (even non symmetric), and therefore zero  since polynomials are dense for the Gaussian measure.

Since $\Gamma_\bE(a_i,a_j)$ are bilinear functions of the $(a_i)$, and $\LL_{\mathbb{OU}}(a_i)$ are linear,  it is clear from the change of variable formula that, for any $k\in \bN$, if $Q(a_i)$ is a polynomial of total degree less than $k$, then so is $\LL_{\mathbb{OU}}((Q)$. Since the set of polynomials in the variables $a_i$ are dense in $\cL^2(\mu)$,  the operator $\LL_{\mathbb{OU}}$ may be then diagonalized  in a basis formed of orthogonal polynomials in the variables $(a_0, \cdots, a_{n-1})$.

\brmq As a consequence, it is worth to observe that if $Q$ is the mean value of $P$, i.e. $Q:=\lag P\rag= \int P(X) d\mu(P)$, where $\mu$ is the image of the gaussian measure, then since $\int \LL_{\mathbb{OU}}((P) d\mu(P)=0$
$$XQ'= nQ.$$Then, $Q(X)= X^n$, which was obvious from the explicit expression of the $a_i$ in terms of $x_i$, which are independent centered Gaussian. But we shall see later  (Remark~\ref{rmk.mean.OU}) that the same computation performed on symmetric matrices leads to a more interesting  result.
\ermq
\section{The image of the spherical Laplacian under elementary symmetric functions\label{sec.diag.sph}}

In an Euclidean $n$-dimensional space, the spherical Laplace   operator $\Delta_{\bS}$ may be written as  the restriction to the unit sphere $\|x\|^2=1$ of some combinations of $\Delta_{\bE}$ and  the Euler operator $D(f) = \sum_i x_i \partial_i f$, namely 
$$\Delta_{\bS} (f)= \Delta_{\bE}(f)-D^2(f)-(n-2)Df.$$
Indeed, when considering the restriction to the sphere of coordinate $x_i$ as a function $\bS^{n-1}\mapsto \bR$, one may describe the spherical Laplace operator $\Delta_\bS$ from
$$\Gamma_\bS(x_i,x_j)= \delta_{ij} -x_ix_j, ~\Delta_\bS(x_i)= -(n-1)x_i.$$
From this, when considering a smooth function $F(x_1, \cdots , x_n)$, the change of variable formula gives  immediately
$\Delta_\bS(F) = \Delta_{\bE}(f)-D^2f-(n-2)Df$.

In the study of Ornstein-Uhlenbeck processes, we already computed the action of $D$ on $P(X)$.  
Then, we get

\bprop  For the Laplace-Beltrami operator $\Delta_\bS$ acting on the unit sphere $\bS^{n-1} \subset \bR^n$, and for the polynomial 
$P(X)= \prod(X-x_i)$, one has
\begin{enumerate}
\item\label{spherdrift}
$\Gamma_\bS(\log P(X),\log P(Y))=\displaystyle{ \frac{1}{Y-X}\Big(\frac{P'(X)}{P(X)}-\frac{P'(Y)}{P(Y)}\Big)- \Big(n-X\frac{P'(X)}{P(X)}\Big)\Big(n-Y\frac{P'(Y)}{P(Y)}\Big)}$
\item\label{sphergamma}
$\Delta_\bS(P(X))= -2n(n-1)P(X)+ 3(n-1)XP'(X)-X^2P''(X).$
\end{enumerate}
\eprop
\bpf
From 
$$\Gamma_\bS(P(X),P(Y))= \Gamma_{\bE}(P(X),P(Y))-DP(X)DP(Y),$$ 
we get point \ref{spherdrift}. Moreover,
$$\Delta_\bS(P(X)= -D(nP(X)-XP'(X))= -(n(nP(X)-XP'(X))+ D(\sum_i ia_iX^i).$$
And
$$D(\sum_i i a_iX^i)= \sum_i i(n-i)a_i X^i= nXP'(X)-X^2P''(X)-XP'(X)\,,$$
which gives point \ref{sphergamma}.

\epf

The support of the image measure is included in the set $\{a_{n-1}^2-2a_{n-2}=1\}$. It does not have a density with respect to he Lebesgue measure $da_0\cdots da_{n-1}$. To deal with it, it is more convenient to look at a more general operator, defined on the unit ball $\bB\subset\bR^n$, and defined, for any $p>1$ from 
$$\Gamma_\bB(x_i,x_j)= \delta_{ij}-x_ix_j, \quad \Delta (x_i)= -px_i.$$ 
Indeed, for any integer $p>n-1$, this operator is nothing else  than the projection on the unit disk of the previous spherical operator on the sphere $\bS^{p}$.
If one observes that
$$\Gamma_\bB(1-\|x\|^2, x_i)= -2(1-\|x\|^2)x_i,$$ comparing with equation~\eqref{eq.rho},  which gives for the density $\rho$ of the associated invariant measure
$\Gamma(\log \rho, x_i)= (n+1-p)x_i$, one sees that this operator has invariant measure 
$$C_{p,n}(1-\|x\|^2)^{(p-n-1)/2} dx_1\cdots dx_n$$ on the unit ball $\bB_n$.

Then, letting $p$ converging to $n-1$, the measure concentrates on the uniform measure on the unit sphere $\bS^{n-1}$. 

One may then consider the image measure through $(a_0, \cdots, a_{n-1})$ of this last operator, which is written as
$$\bcas\Gamma_\bB(\log P(X),\log P(Y))= \frac{1}{Y-X}\Big(\frac{P'(X)}{P(X)}-\frac{P'(Y)}{P(Y)}\Big)- \Big(n-X\frac{P'(X)}{P(X)}\Big)\Big(n-Y\frac{P'(Y)}{P(Y)}\Big),\\
\\
 \Delta_\bS(P(X))= -n(p+n-1)P(X)+ (p+2(n-1))XP'(X)-X^2P''(X).\ecas$$
 Its reversible  measure is then, up to some scaling factor 
 $$\I_{\cD} D(a_0,\cdots,  a_{n-1})^{-1/2} (1+2a_{n-2}-a_{n-1}^2)^{(n-p-1)/2}da_0\cdots da_{n-1},$$ which  concentrates on the set $1+2a_{n-2}-a_{n-1}^2=0$ when $p\to n-1$.
 
 \section{ Symmetric  operators on the spectrum of real symmetric matrices\label{sec.sym.matr}}
 
 The space of symmetric real matrices is an Euclidean space with the norm $\|M\|^2= \tr (M^2)$. The spectrum of such matrices  is real and we want to describe the action of the Euclidean  and spherical Laplacians  and the Ornstein-Uhlenbeck  on the spectrum of the matrix. For the moment, we only deal with the Euclidean Laplacian.
 
 For this, we use first start with the description of the Euclidean Laplace operator on the entries $M_{ij}$ of the matrix.
 
 One has 
 $$\Gamma_{\bE,\bR}(M_{ij}, M_{kl})= \frac{1}{2} (\delta_{ik}\delta_{jl}+ \delta_{il}\delta_{jk}), ~\LL_{\bE,\bR}(M_{ij})=0.$$ 
 This formula just captures the fact that the entries of the matrix are independent Brownian motions, subject to the restriction that the matrix is symmetric. 
More precisely, it   encodes by itself the fact the associated process lives in the space of symmetric matrices. Since we shall use this kind of argument in many places,  it is worth to explain why in this simple case.  In terms of the  stochastic process $(X_t)$ associated with $\LL_{\bE,\bR}$, which lives a priori in the space of matrices,  and for the functions $h_{ij}(x)= M_{ij}-M_{ji}$, it is not hard to see that $\LL_{\bE,\bR}(h_{ij}) = 0$ and  $\Gamma(h_{ij}, h_{ij})=0$. Therefore, for the associated Brownian motion $(X_t)$,  and any smooth function $f$, $\LL_{\bE,\bR}(f(h_{ij}))=0$ and, provided $f(h_{ij})$ and $\LL_{\bE,\bR}(f(h_{ij}))$ are bounded,  $\E(f(h_{ij}(X_t) \slash X_0=x)= f(h_{ij})(x)$. This shows that $h_{ij}(X_t)$ remains constant (almost surely).  Then, if the process  starts from a symmetric matrix, it stays forever in the space of symmetric matrices (we shall not need this fact for what follows).
 
 We start with the following elementary lemma, which will be in full use through the rest of the paper
 
 \blem Let $M= (M_{ij})$ be a matrix and $M^{-1}$ be its inverse, defined on the set $\{\det M\}\neq 0$. Then
 \benum
 \item \label{deriv.log.det}$\partial_{M_{ij}} \log \det M= M^{-1}_{ji}$, ~
 \item \label{deriv.inverse.mat2} $\partial_{M_{ij}}\partial_{M_{kl}}\log \det M= -M_{jk}^{-1}M_{li}^{-1}$.
 \eenum
 \elem
 \bpf
 We first start with 
 \beq\label{deriv.inverse.mat}\partial_{M_{ij}} M^{-1}_{kl}= -M_{ik}^{-1}M_{lj}^{-1}.\eeq
 To see this, consider  the formula
 $\sum_k M^{-1}_{ik}M_{kj}= \delta_{ip}$, that we derive with respect to $M_{pq}$, leading to
$$\sum_k (\partial_{M_{pq}}M^{-1}_{ik})M_{kj}+ M^{-1}_{ik}\delta_{kp}\delta_{jq}=0.$$
Fixing $p$ and $q$, if $DM$ denotes the matrix $\partial_{M_{pq}}M^{-1}$, one gets for every $p,q$ 
$(DM\times M)_{ij}+ M^{-1}_{ip}\delta_{jq}=0$, which we may now multiply from the right by $M^{-1}$ to get~\eqref{deriv.inverse.mat}.

Now, we observe that $P_{ij}= \det(M)M^{-1}$ is the comatrix,  and therefore the $M_{ji}$ entry of this matrix is a polynomial in the entries of the matrix which does does not depend on $M_{ij}$.

One gets
$$(\partial_{M_{ij}} \det M)  M^{-1}_{ji}- \det (M) M^{-1}_{ji}M^{-1}_{ji}=0,$$ which is an identity between rational functions in the entries $M_{ij}$. Therefore, we may as well divide both terms by $M_{ij}^{-1}$ to obtain item~\ref{deriv.log.det}. Item~\ref{deriv.inverse.mat2} follows immediately.

 \epf
 We are now in position to consider the action of the Laplace operator on the spectral measure of $M$. If we consider  the characteristic polynomial $P(X)= \det(M-X\Id)$, one has
 \bprop
 \begin{enumerate} 
 \item
$  \Gamma_{\bE,\bR}(\log P(X), \log P(Y))= \frac{1}{Y-X}\Big(\frac{P'(X)}{P(X)}-\frac{P'(Y)}{P(Y)}\Big)$
\item
$\LL_{\bE,\bR} P(X)= -\frac{1}{2}P''.$
 \end{enumerate}
 In other words, the diffusion associated to the spectrum of $M$ has the same operator carré du champ, and the invariant measure has a constant density with respect to the Lebesgue measure $\I_{D}da_0\cdots da_{n-1}$.
 
 \eprop
 
 \bpf We start with the first formula.
 Fix $X$ and $Y$  and write $U= (M-X\Id)^{-1}$, $V= (M-Y\Id)^{-1}$.    One has
 $$\Gamma_{\bE,\bR}(\log P(X), \log P(Y))= \sum_{ijkl} \partial_{M_{ij}} \log P(X) \partial_{M_{kl}}\log P(Y) \Gamma_{\bE,\bR}(M_{ij}, M_{kl}).$$
 Since $P(X)= \det (M-X\Id)$, and from the value of $\Gamma_{\bE,\bR}(M_{ij},M_{kl})$, this writes
 $$\frac{1}{2}\sum_{ijkl} U_{ji}V_{lk}(\delta_{ik}\delta_{jl}+ \delta_{il}\delta_{jk})= \frac{1}{2}\sum_{ijkl} (U_{ji}V_{ji}+U_{ji}V_{ij})= \tr (UV).$$
 But, if $x_i$ are the eigenvalues of $M$, then $P(X)= \prod_i (x_i-X)$, $P(Y)= \prod_i(x_i-Y)$, and 
 \beqnas \tr (UV) &=& \sum_i \frac{1}{(x_i-X)(x_i-Y)}= \frac{1}{X-Y} \sum_i \frac{1}{x_i-X} -\frac{1}{x_i-Y}\\ &=& \frac{1}{Y-X} \Big(\frac{P'(X)}{P(X)}-\frac{P'(Y)}{P(Y)}\Big).\eeqnas
 
 It remains to compute 
 $$\LL_{\bE,\bR}(\log P(X))= \sum_{ij} \partial_{M_{ij}} \log P(X) \LL_{\bE,\bR}(M_{ij})+ \sum_{ijkl} \partial_{M_{ij}}\partial_{M_{kl}} \Gamma_{\bE,\bR}(M_{ij}, M_{kl}),$$ which writes
 $$-\frac{1}{2}\sum_{ijkl}U_{jk}U_{li}(\delta_{ik}\delta_{jl}+\delta_{il}\delta_{jk})= -\frac{1}{2}\big(\tr (U^2)+ (\tr U)^2\big).$$
 
 Writing 
 $$\frac{\LL_{\bE,\bR} P(X)}{P(X)}= \LL_{\bE,\bR} \log P(X)+ \Gamma_{\bE,\bR}(\log P(X), \log P(X)),$$ and noticing that
 $$\tr (U^2)= \frac{P'(X)^2}{P(X)^2}- \frac{P''(X)}{P(X)}= \Gamma_{\bE,\bR}(\log P(X), \log P(X)), ~\tr U= -\frac{P'(X)}{P(X)}$$
 one gets the formula for $\LL_{\bE,\bR} P(X)$.
 
 Comparing with  Corollary \ref{cor.gamma.discr}, one sees that $\LL_{\bE,\bR}(P)= \frac{1}{2} \Gamma_{\bE,\bR}(P, \log \discr(P))$. The last point then  is just the observation that if $\rho$ is the reversible measure for the image of $\LL_{\bE,\bR}$, then $\rho$ has, up to a constant,  density $\discr(P)^{1/2}$ with respect to the Riemann measure, which is just $C\discr(P)^{-1/2}$.
 
 \epf
 
 \brmq Moving back the measure to the Weyl chamber $\{\lambda_1<\cdots< \lambda_n\}$, one sees that the density of the spectral measure with respect to $d\lambda_1\cdots d\lambda_n$ is $C\discr(P)^{1/2}= C\prod |\lambda_i-\lambda_j|$.
 
 \ermq
 
 If one wants to extend the previous computation to the Gaussian or spherical case, one has to consider also the image of $D= \sum_{ij} M_{ij} \partial _{M_{ij}} =\frac{1}{2}\Gamma_{\bE,\bR}(\|M\|^2,\cdot)$ on the spectral function $P(X)= \det(M-X\Id)$. 
 One gets, with $U(X)= (M-X\Id)^{-1}$, 
 \beqnas D(\log P(X))&=& \sum_{ij} M_{ij} \partial_{M_{ij} } P(X)= \sum_{ij}M_{ij} U(X)_{ji}= \tr MU(X)\\&=& \tr (\Id+ XU(X))= n -X\frac{P'}{P},\eeqnas
 from which $D(P) = nP- XP'$.
 Therefore, if 
 $P(X) = \sum_i a_iX^i$, $D(a_iX^i)= (n-i) a_i$. If one wants to consider also the action of the spherical Laplace operator, one needs also to consider $D^2(P)$.
 But $D^2(a_iX^i)= (n-i)^2a_i$, from which 
 $$D^2(P)= n^2P -(2n-1)XP' + X^2P'',$$
 that is with $DP= (n\Id-X\partial_X )P$, $D^2(P)= (n\Id-X\partial_X)^2 P$, although there is no reason a priori for this last identity, since $DP$ is no longer the characteristic polynomial of the symmetric matrices.  Observe that the action of $D$ on $P$ is similar that the one in the diagonal case.  This is not surprising since setting $U= a_{n-1}^2-2a_{n-2} = \|M\|^2$, 
 $D(F)= \frac{1}{2} \Gamma_{\bE,\bR}(U,F)$, for any function $F$.
 
 Then, one has for the spherical operator on symmetric matrices, $\LL_{\bS,\bR}= \LL_{\bE,\bR}- D^2-(N-2)D$, where $N= n(n+1)/2$, from which
 $$\LL_{\bS,\bR}(P)= -(\frac{1}{2}+X^2) P'' + \frac{(n+6)(n-1)}{2} XP' -\frac{n(n+4)(n-1)}{2} P.$$ and 
 $$\Gamma_{\bS,\bR}(P(X), P(Y))= \frac{1}{Y-X} (P'(X)P(Y)-P'(Y)P(X)) - (nP(X)-XP'(X))(nP(Y)-XP'(Y)).$$
 
 \brmq\label{rmk.mean.OU}
  If we perform the same computation from the Gaussian measure $\gamma$ instead of the  Lebesgue one  (that is if we start from  an Ornstein-Uhlenbeck process instead of the Brownian motion on matrices), we end up with 
  $\LL_{\mathbb {OU},\bR}(P)= -\frac{1}{2} P'' + XP'-nP$.
  
  Now, if  we consider now  the mean value polynomial, that is $Q= \lag P\rag= \int Q d\gamma$, one gets 
  $Q''-2XQ' = -2nQ$. From which we see that, up to some constant, $Q(X/\sqrt{2})$ is an Hermite polynomial.
  
  On the other hand,  the same computation for the spherical case leads to
  $$(\frac{1}{2} +X^2) Q''  -\frac{(n+6)(n-1)}{2} XQ' = -\frac{n(n+4)(n-1)}{2} Q. $$ 
  In the same way that  Hermite polynomials are  the orthogonal polynomial family   associated with Gaussian measure, one would expect some connection between those polynomials $Q$ and the one-dimensional projection of the uniform  measure on a sphere (in some dimension), i.e. Jacobi polynomials, but it does not seem to be the case.
 
 \ermq

  \section{  Symmetric  operators on the spectrum of Hermitian  and quaternionic   matrices\label{sec.herm.matr}}
  
  \subsection{Hermitian matrices}
  
  In this section, we extend the previous computations to Hermitian matrices. We mainly consider an Hermitian matrix on $\bC^n$ as a real symmetric matrix on $\bR^{2n}$. Indeed, considering a vector $Z= S+iT$ in $\bC^n$,  where $S$ and $T$ are the real and imaginary part of $Z$, an Hermitian matrix  $H$ may be  seen as  $M+ iA$, where $M$ is an $n\times n$ real symmetric and $A$ is $n\times n$ real antisymmetric matrix.  Then, writing $H$ as a bloc matrix, we have
  $H= \bpm M& A\\-A&M\epm$, which is a real  symmetric matrix with special structure. Indeed, any real eigenspace for this matrix is at least 2 dimensional (and is exactly  2-dimensional in the generic case), since if $Z= (S,T)$ is an eigenvector, so is $(-T,S)$, which corresponds to the eigenvector $iZ$.
  
  Moreover, the determinant $P(X)$  of $H-X\Id$  may be written as $Q(X)^2$, where $Q(X)$ is  a polynomial whose coefficients are polynomials in the entries of $M$ and $A$ (actually, $Q(X)$ is the Pfaffian of the anti-symmetric $2n\times 2n$ matrix $\bpm A &M-X\Id\\-M+X\Id &A\epm$. Therefore, if we consider the entries of $M$ and $A$ uniformly distributed under the Lebesgue measure, as we did for real symmetric ones,  the spectrum of the matrix $H$ is certainly not absolutely continuous with respect to the Lebesgue measure, and we are more interested indeed in the law of the roots of $Q$ than in the law of the roots of $P$. 
  
  As before, we look at the Euclidean Laplace  operator $\LL_{\bE,\bC}$ acting on $H= M+iA$, with $M=(M_{ij})$ and $A= (A_{ij})$ and we encode the symmetries via the following formulae
  $$\bcas \LL_{\bE,\bC}(M_{ij})= \LL_{\bE,\bC}(A_{ij})=0\\
  \Gamma_{\bE,\bC}(M_{ij}, M_{kl})= \frac{1}{2} (\delta_{ij}\delta_{kl}+ \delta_{ik}\delta_{jl})\\
  \Gamma_{\bE,\bC}(A_{ij}, A_{kl})= \frac{1}{2} (\delta_{ij}\delta_{kl}-\delta_{ik}\delta_{jl})\\
  \Gamma_{\bE,\bC}(M_{ij}, A_{kl})= 0
  \ecas
  $$  
 It is worth to observe that any power (and therefore the inverse when it exists) of an Hermitian matrix is again an Hermitian matrix, and one may perform the same computation as before on $P(x)= \det (H-X\Id)$ (still considered as a $2n\times 2n$ matrix as before).
 
 We obtain, with again $U(X)= (H-X\Id)^{-1}$,

 \bprop\label{gen.herm.2}
 
 $$\Gamma_{\bE,\bC}(\log P(X), \log P(Y))= 2 \tr(U(X)U(Y))= \frac{2}{Y-X}(\frac{P'(X)}{P(X)}-\frac{P'(Y)}{P(Y)}).$$
 
 $$\LL_{\bE,\bC}(P(X))= \frac{3}{2}\frac{P'(X)^2}{P(X)}-2P''(X).$$

 \eprop 
 
 We do not give the details of the proof here since we shall give a more general result in the setting of Clifford algebras, for which this is just the simplest  example.
 
 But is is worth to observe the following. Since 
 $P(X)= \sum_i X^i a_i$, $\LL_{\bE,\bC}(P(X))= \sum_i X^i\LL_{\bE,\bC}(a_i)$.
 Therefore, this last quantity has to be a polynomial, and then  $\frac{P'(X)^2}{P(X)}$ is a polynomial in $X$. This implies in particular that all the roots of $P$ have multiplicities at least $2$, since every root of $P$ is also a root of $P'$. In particular, the image measure of the Lebesgue measure is not absolutely continuous with respect to $d\mu_0$.
 
Observe furthermore that  if we set 
 $P= P_1^2$, then one gets from the change of variable formula
 \beq\label{gen.herm.1}\Gamma_{\bE,\bC}(\log P_1(X), \log P_1(Y))= \frac{1}{Y-X} \Big( \frac{P_1'(X)}{P_1(X)}-\frac{P_1'(Y)}{P_1(Y)}\Big),~\LL_{\bE,\bC}(P_1)= -P_1''.\eeq
 
 In particular, moving back to the Weyl chamber $\{ \lambda_1 <\cdots < \lambda_n\}$, the invariant measure is, up to a constant, $\prod_{i<j} (\lambda_i-\lambda_j)^2 d\lambda_1 \cdots d\lambda_n$.
 
 Are we able to deduce directly from the form of the generator that indeed $P(X)$ has almost surely double roots? We shall see that it is indeed  the case. It may be seen at this level as a purely formal argument, since we know in advance that in this  Hermitian case the roots are double.  But later we shall face similar situations, where we do not know in advance the multiplicity of the roots, and we want to be able to deduce them from the generator.    More precisely,  we shall see that if a generator of the form  given in Proposition~\ref{gen.herm.2} maps polynomials into polynomials, then those polynomials must have roots with multiplicity 2.  This relies on Lemma~\ref{lem.puiss.P} and Proposition~\ref{lem.puiss.P2}.

From the form of the operator, one already sees that there are some algebraic relations between the coefficients $a_i$ of the polynomial $P$. The following Lemma~\ref{lem.puiss.P}  is quite formal, and allows to devise from the form of the operator the multiplicity of the roots of $P$. 
Then, Proposition~\ref{lem.puiss.P2} provides a proof that the multiplicity of the roots are   indeed what is expected.

 \blem\label{lem.puiss.P}  Suppose that a  diffusion operator  $\LL$ on some set of analytic functions $P(X)= \sum_i a_i X^i$ in the variable $X$  satisfies, for some constants $\alpha, \beta, \gamma$,
 \beq\label{eq.gal.P}\LL(P)= \alpha P''+ \beta \frac{P'^2}{P}, ~\Gamma(\log P(X), \log P(Y))=  \frac{\gamma}{Y-X} \Big( \frac{P'(X)}{P(X)}-\frac{P'(Y)}{P(Y)}\Big).\eeq Let  $a\in \bR$, $a\neq 0$, and set $P= P_1^a$. 
 Then
 $$\Gamma(\log P_1(X), \log P_1(Y))= \frac{\gamma/a}{Y-X} \Big( \frac{P_1'(X)}{P_1(X)}-\frac{P_1'(Y)}{P_1(Y)}\Big)$$
and 
  $$\LL P_1= (\alpha+\gamma\frac{a-1}{a}) P_1'' + (a(\alpha+ \beta)+ \gamma \frac{1-a}{a}-\alpha) \frac{P_1'^2}{P_1}.$$
In particular, if $a$ satisfies 
 \beq \label{eq.puiss}a^2(\alpha+\beta)-a(\alpha+ \gamma) + \gamma=0,\eeq
 then, 
 $$\LL(P_1)= a(\alpha+ \beta) P_1''.$$

 \elem 
 
 Therefore, one may expect that, whenever $\LL$ maps polynomials into polynomials, and  precisely for those  values of $a$ solutions of equation~\eqref{eq.puiss}, the roots of $P$ have multiplicity $a$.
  
 \bpf The formula for $\Gamma(\log P_1(X), \log P_1(Y))$ is immediate. The formula for $\LL(P)$ 
  follows easily from the remark that
 $$\frac{\LL P}{P}= \alpha\partial^2_X \log P+ (\alpha+\beta)(\partial_X \log P)^2, \Gamma(\log P, \log P)= -\gamma \partial_X^2 \log P.$$
 Then, 
 $$a\frac{\LL P_1}{P_1}= \frac{\LL P}{P} + \frac{1-a}{a}\Gamma(\log P, \log P),$$ and this gives
 $$\frac{\LL P_1}{P_1}= (\alpha+\gamma\frac{a-1}{a}) \frac{P_1''}{P_1} + (a(\alpha+ \beta)+ \gamma \frac{1-a}{a}-\alpha) \frac{P_1'^2}{P_1^2}.$$

 \epf
In the next proposition~\ref{lem.puiss.P2}, we consider polynomials $P(X)$ with coefficients $a_i$ which are polynomials in some variables $(x_i)$ (in our case the entries of a matrix). Then, when writing $P(X)= \prod(X-\lambda_i)^{\alpha_i}$, with $\lambda_1<\lambda_2< \cdots <\lambda_k$, the multiplicities $\alpha_i$ may only change on some algebraic surface in the set of coefficients $(x_i)$. Those algebraic surfaces having Lebesgue measure $0$, and our operators $\LL$ being local, we may as well (up to some localization procedure and outside a set of Lebesgue measure $0$) consider them as constants.

 \bprop\label{lem.puiss.P2} Let $\LL$ be a diffusion  operator acting on a set  of degree $d$ monic   polynomials, with values in the set of degree $d$ polynomials,  and satisfying equation~\eqref{eq.gal.P}.  Then, every root of $P$ has multiplicity $\alpha_1$ or $\alpha_2$, where $\alpha_i$,  $i=1,2$ are the roots of equation~\eqref{eq.puiss}.
 
 In particular,equation~\eqref{eq.gal.P} may  only hold for polynomials whenever equation~\eqref{eq.puiss} has at least one integer solution.

 \eprop
 
 In practise, equation~\eqref{eq.puiss} will have only one integer positive  root, which will allow us to identify the multiplicity without any ambiguity. Moreover, in this situation, we may set $P= P_1^a$, where $P_1$ is a polynomial,  and Lemma~\ref{lem.puiss.P} applies within the set of polynomials.
 
 \bpf Let us consider $\lambda_1<\lambda_1\cdots < \lambda_k$ the different roots of $P(X)$, and set 
 $$P(X) = \prod_1^k (X-\lambda_i)^{\alpha_i},$$ where $\alpha_i\geq 1$.
 Then
 $$\Gamma(\log P(X), \log P(Y))= \sum_{ij} \frac{\alpha_i\alpha_j}{(X-\lambda_i)(Y-\lambda_j)} \Gamma(\lambda_i, \lambda_j).$$
 On the other hand, we know  from Lemma~\ref{lem.puiss.P} that 
 $$\Gamma(\log P(X), \log P(Y))=  \frac{\gamma}{Y-X}  \Big( \frac{P'(X)}{P(X)}-\frac{P'(Y)}{P(Y)}\Big),$$
 which translates into
 $$\Gamma(\log P(X), \log P(Y))= \frac{\gamma}{Y-X}\sum_i \frac{\alpha_i}{X-\lambda_i}-\frac{\alpha_i}{Y-\lambda_i}.$$
 Identifying both expressions leads to
 $$\Gamma(\lambda_i, \lambda_j)= \gamma\frac{\delta_{ij}}{\alpha_i}.$$
  
 Also, on the one hand,
 $$\LL(\log P)= \sum_i \frac{\alpha_i}{X-\lambda_i} \LL \lambda_i -\sum_i \frac{\alpha_i}{(X-\lambda_i)^2} \Gamma(\lambda_i, \lambda_i),$$ and on the other, from~\eqref{eq.gal.P}
 \beqnas \LL(\log P)&=& \frac{ \LL P}{P} -\Gamma(\log P, \log P)\\
 &=& (\alpha+\gamma) \partial_X^2\log P+ (\alpha+\beta)(\partial_X \log P)^2\\
 &=&-(\alpha+\gamma) \sum_i \frac{\alpha_i}{(X-\lambda_i)^2}+  (\beta+\alpha) (\sum_i \frac{\alpha_i}{X-\lambda_i})^2.\eeqnas
 Identifying the terms in $(X-\lambda_i)^{-2}$ leads to
 $$\alpha_i^2(\beta+\alpha) -(\alpha+\gamma) \alpha_i+ \gamma=0.$$

 \epf

 In particular, applying Lemma~\ref{lem.puiss.P} in the case of Hermitian matrices leads, with $\alpha= -2, \beta= 3/2$ and $\gamma=2$ to 
 $\alpha_i=2$, and then, setting $P= P_1^2$, to~\eqref{gen.herm.1} for $\LL_{\bE,\bC}(P_1)$ and $\Gamma_{\bE,\bC}(\log P_1(X), \log P_1(Y))$. 
This  in turns shows that every root of $P$ has multiplicity $2$, and  therefore that $P$, up to some sign may be written $P_1^2$, where $P_1$ is a polynomial for which~\eqref{gen.herm.1} holds.
  As a consequence,  the image measure for the roots of $P_1$ has density $\discr(P_1)$ wit respect to $d\mu_0$.

\brmq It is  worth to observe that if $P$ is a monic polynomial whose coefficients are polynomials in some variables $(x_1, \cdots, x_n)$, and if $P= P_1^a$, where $P_1$ is a polynomial, then  $P_1$ is monic and the coefficients of $P_1$  are again polynomials in the variables $(x_1, \cdots, x_n)$. In the case of Hermitian matrices, this shows in particular that the determinant of a matrix of the form $\bpm M&A\\-A&M\epm$, where $M$ is symmetric and $A$ is antisymmetric may be written as $Q^2$, where $Q$ is a polynomial in the entries of $M$ and $A$ (indeed, it is nothing else up to a sign than the Pfaffian of $\bpm A &M\\-M&A\epm$, but it is worth deriving it by  pure probabilistic arguments). 

\ermq
\subsection{Quaternionic matrices}

Here, we are given $M,A^1,A^2,A^3$ where $M$ is symmetric and $A^i$ are antisymetric.

The associated real symmetric matrix is then
$$\cM=\bpm M&A^1&A^2&A^3\\
-A^1&M&A^3&-A^2\\
-A^2&-A^3&M&A^1\\
-A^3&A^2&-A^1&M\epm$$
The eigenspaces are 4-dimensional and the determinant of such a matrix may be written $Q^4$, where $Q$ is a polynomial in the entries of the various matrices.

A real $4n\times 4n$ matrix $\cM$ having this structure will be called a $\cH$-symmetric matrix. It is quite immediate that if $\cM$ is $\cH$-symmetric, such is $\cM^k$ for any $k\in \bN$, and also such is $\cM^{-1}$ on the set where $\cM$ is invertible.

On the entries of $M_{ij}$ and $A^k_{ij}$, we shall impose the metric coming from the euclidean metric on $\cM$. This gives
$$\bcas \Gamma_{\bE,\bH}(M_{ij}, M_{kl})= \frac{1}{2} (\delta_{ij}\delta_{kl}+ \delta_{ik}\delta_{jl})\\
\Gamma_{\bE,\bH}(A^p_{ij},A^q_{kl}) = \frac{\delta_{pq}}{2}((\delta_{ij}\delta_{kl}- \delta_{ik}\delta_{jl})\\
\Gamma_{\bE,\bH}(M_{ij}, A^p_{kl})=0\ecas.$$
We also impose
$$\LL_{\bE,\bH} M_{ij}= \LL_{\bE,\bH}A^p_{ij}=0.$$

Setting $U(X)= (\cM-X\Id)^{-1}$, $P(X)= \det(\cM-X\Id)$, one has 
 $$\Gamma_{\bE,\bH}\big(P(X),P(Y)\big)= \frac{4}{Y-X}\big(P'(X)P(Y)-P'(Y)P(X)\big),$$
$$\frac{ \LL_{\bE,\bH}(P)}{P}= \frac{\Gamma_{\bE,\bH}(P,P)}{P^2}+\tr(U(X)^2) - \frac{1}{2} (\tr U(X))^2.$$

And in the end
$$\LL_{\bE,\bH}P= \frac{9}{2}\frac{P'^2}{P}-5 P''.$$

Looking for which $a$, one has $P= P_1^a$, equation~\eqref{eq.puiss} on $a$ leads to
$a^2-2a-8=0$, for which the unique positive solution is $a=4$, leading to
$$\LL_{\bE,\bH} P_1= -2P_1'', ~\Gamma_{\bE,\bH} (\log P_1(X), \log P_1(Y))= \frac{1}{Y-X}\Big( \frac{P_1'(X)}{P_1(X)}-\frac{P_1'(Y)}{P_1(Y)}\Big),$$
 and Proposition~\ref{lem.puiss.P2} shows that  all the roots of $P$ have multiplicity $4$, and that $P_1$ is indeed a polynomial. 
 
 In the end, one obtains that the reversible measure for the image operator has density  is $\discr(P)^2$ with respect to the Riemann measure, or in other terms $\discr(P)^{3/2}$ with respect to the measure $da_0\cdots da_{n-1}$. Back to the Weyl chamber $\{\lambda_1< \lambda_n\}$, the invariant measure has density $\prod_{i<j} (\lambda_i-\lambda_i)^4$ with respect to the Lebesgue measure.
 
 \section{Symmetric matrices on general Clifford algebras\label{sec.gal.cliff}}
 There are many natural algebras with dimension $2^p$. Among them, let us mention  exterior algebras, Cayley-Dickson algebras and Clifford algebras. Since $2^p$ is the cardinal of $\cP(\{1, \cdots, p\})$, it is natural to look for a basis $\omega_A$ for such algebras, where $A\subset E$, and $|E|= n$. If we denote by $A\DS B$ the symmetric difference $A\cup B \setminus (A\cap B)$, in those three cases one has
 $\omega_A \omega_B = (A|B) \omega_{A\DS B}$, where $(A|B)$ takes values in $\{-1,0,1\}$.
 
We define general Clifford 	algebras  are the ones where the algebra is associative and $(A|B)\in\{-1,1\}$.  We shall impose $\omega_\emptyset$ to be the unitary element  of the algebra. The associativity imposes that, for any triple $(A,B,C)$ of elements of $\cP(E)$, one has
$$(A|\emptyset)=(\emptyset|A)=1, ~(A|B\DS C)(B|C)= (A|B)(A\DS B|C).$$

It is worth to reduce to the case $E= \{1, \cdots, p\}$ (that is to decide that $E$ is an ordered set), such that up to a change of sign in $\omega_A$, one may always suppose that $\omega_A= \omega_{i_1}\cdots \omega_{i_k}$ when $A= (i_1, \cdots , i_k)$. Therefore, one sees that all the multiplication rules are just given by $e_ie_j= \pm e_je_i$ and $e_i^2= \pm e_i$. In which case, we are reduced to $(i|j)= 1$ if $i<j$ and $(A|B)= \prod_{i\in A, j\in B} (i|j)$, from which we get
$$(A\DS B|C)= (A|C)(B|C), ~(A|B\DS C)= (A|B)(A|C).$$

The general Clifford algebra is then just determined by the choice of the various signs in $(i|j)(j|i)$ for $i<j$ and $(i|i)$. But many such different choices may give rise to isomorphic algebra : for example, given any Clifford algebra and any choice $(A_1, \cdots, A_p)$ which generates $\cP(E)$ by symmetric difference would produce a Clifford algebra isomorphic to the starting one with signs $(A_i|A_j)$ instead of $(i|j)$ (think for example of $A_i= \{1, \cdots, i\}$).

The Clifford algebra $\Cl(E)$ is then $\{\sum_{A\subset E} x_A \omega_A, x_A\in \bR\}$ that we endow with the standard Euclidean metric in $\bR^{2^p}$(that is $(\omega_A, A\in \cP(E))$ form an orthonormal basis).

Now, we consider on $\bR^n\otimes\Cl(E)$ matrices $\sum_A M^A\omega_A$, where $M^A$ are $n\times n$ matrices, acting on $\bR^n\otimes \Cl(E)$ by 
$$(\sum_A M^A\omega_A)(\sum_B X^B\omega_B)= \sum_{A,B} M^A X^B(A|B)\omega_{A\DS B}= \sum_{A,B} (A\DS B|B)M^{A\DS B}X^B\omega_A,$$ and we end up with bloc matrices $(M^{A,B}_{ij})$, where 
$M^{A,B}= (A\DS B|B) M^{A\DS B}$.  Indeed, what we did is to associate to a matrix $M$ with coefficient in the algebra $\Cl(E)$ a matrix $\phi(M)$ with real coefficients, in a linear injective way. It turns out that,  thanks to the associativity of the the algebra $\Cl(E)$, this is an algebra homomorphism, that is $\phi(MN)= \phi(M)\phi(N)$. 

 Endowing $\bR^n\otimes \Cl(E)$ with the associated Euclidean metric, we may therefore look at those matrices $\phi(M)$  which are symmetric,  one sees that the requirement is that 
 $(M^{A})^t= (A|A) M^A$, and then the associated bloc matrix is 
 $\cM=((A\DS B|B)M^{A\DS B})$. We shall call those symmetric matrices $\Cl(E)$-symmetric matrices.

 We now chose the Euclidean metric on those $\Cl(E)$-symmetric matrices, and look at the associated Laplace operator. One then sets
 \beq\label{eq.LGamma.clif.sym}\Gamma_{\bE,\Cl}(M^A_{ij},M^B_{kl})=   \frac{1}{2}\delta_{A,B} (\delta_{ik}\delta_{jl}+ (A|A) \delta_{il}\delta_{jk}), \LL (M^A_{ij})=0.\eeq

Going to the associated stochastic processes, these formulae just say that the various entries of the matrices  are independent Brownian motions, subject to the restrictions that the matrices must satisfy the symmetry relations that we have just  which are imposed by the algebra structure of $\Cl(E)$.

The aim is now to compute when possible the image of this Laplace operator on the spectrum of $\cM$. We shall see that is strongly depends on the sign structure of the algebra $\Cl(E)$. In the next Section,   we shall  reduce our analysis to standard Clifford algebras, that is when $(i|i)= (i|j)(j|i)=-1$, for any $(i,j)\in E^2$. But is is worth to describe  first the  computations in the general case. Indeed,  as mentioned before, since many different sign structures lead to  isomorphic algebras, the various quantities which will appear in the computations will  to be invariant under those isomorphisms, and it is worth to identify them. 

The first task is to observe that, if $\cM$ is $\Cl(E)$-symmetric, so is $\cM^k$ for any $k$. This is a direct consequence of the algebra homomorphism, since if $\cM= \phi(M)$, then $\cM^k= \phi(M^k)$. Therefore, on the set where $\det(\cM)\neq 0$, so is its inverse $\cM^{-1}$ and also $U(X)= (\cM-X\Id)^{-1}$, for $X$ not in the spectrum of $\cM$. Indeed, since for $\|\cM\|$ close to $0$, $(\Id -\cM)^{-1} = \sum_k \cM^k$, and consequently, for $X\neq 0$ and $\cM$ small enough, then, $U(X)$ is $\Cl(E)$-symmetric. Then, since the property of being $\Cl(E)$-symmetric is linear in the coefficients of $\cM$, and since the coefficients of $U(X)$ are rational functions of the coefficients of $\cM$, the $\Cl(E)$-symmetry of $U(X)$ may be extended  from small values of $\cM$ to any $\cM$ which is $\Cl(E)$-symmetric.

Once this is observed,  and still denoting  $U(X)= (\cM-X\Id)^{-1}$, it may be written as a block matrix $((A\DS B|B)U(X)^{A\DS B}_{ij})$, where $U(X)^{A}= U(X)^{A, \emptyset}$ is such that $(U(X)^A)^t= (A|A) U(X)$. Then, the method used for real symmetric matrices  may be extended to $\Cl(E)$-symmetric matrices and we get

\bprop\label{prop.gal.cliff} Let $P(X)= \det(\cM-X\Id)$ and $U(X)= (\cM-X\Id)^{-1}$. Then
$$\Gamma_{\bE,\Cl}(P(X),P(Y))= \frac{2^p}{Y-X}(P'(X)P(Y)-P'(Y)P(X)), $$
and 
$$\frac{\LL_{\bE,\Cl} (P)}{P}= \Gamma_{\bE,\Cl}(\log P)- \frac{1}{2}\Big(\sum_{A\subset E} (A|A)\Big)\tr\big(U(X)^2\big)-2^{p-1} \sum_C (C|C) H(C)\big(\tr U(X)^C\big)^2,$$ where
$$H(C)= \sum_{A\subset E}  (A|C)(C|A).$$

Moreover, 
$$\tr(U(X)^2)= \frac{P'^2}{P^2}-\frac{P''}{P}.$$
\eprop

\bpf
Let us start with the formula for $\Gamma$.
If $U(X)= (U(X)^{A, B}$, where $U(X)^{A,B}= (A\DS B| B) U(X)^{A\DS B}$,  one has, using the change of variable formula and  equation~\eqref{deriv.log.det},  
$$\Gamma_{\bE,\Cl}(\log P(X), \log P(Y))= \sum_{A,B,C,D, i,j,k,l} U(X)^{B,A}_{ji}U(Y)^{D,C}_{lk} \Gamma_{\bE,\Cl}(\cM^{A,B}_{ij}, \cM^{C,D}_{kl}).$$
Now, since 
$\cM^{A,B}= (A\DS B|B) M^{A\DS B}$ and $\cM^{C,D}= (C\DS D|D) M^{C\DS D}$, and from~\eqref{eq.LGamma.clif.sym}, one gets
$$\Gamma_{\bE,\Cl}(\cM^{A,B}_{ij}, \cM^{C,D}_{kl})= I_{\{A\DS B\DS C\DS D= \emptyset\}}(A\DS B|B\DS D)\frac{1}{2}(\delta_{ik}\delta_{jl}+ (A\DS B|A\DS B)\delta_{il}\delta_{jk}).$$
On the other hand
$$\tr U(X)U(Y)= \sum_{A,B,i,j} U(X)^{A,B}_{i,j}U(Y)^{B,A}_{ji}= \sum_{A,B} (A\DS B|A\DS B) \tr U(X)^{A\DS B}U(Y)^{A\DS B}.$$
From this, we get
\beqnas\Gamma_{\bE}(\log P(X), \log P(Y))&=&\sum_{A\DS B\DS C\DS D= \emptyset}
(A\DS B|A\DS B) \tr\Big(U(X)^{A\DS B}U(Y)^{A\DS B}\Big)\\
&=& 2^p\sum_{A,B}(A\DS B|A\DS B)\tr\Big(U(X)^{A\DS B}U(Y)^{A\DS B}\Big)\\&=&
2^p\tr U(X)U(Y) .\eeqnas
If we denote by  $\lambda_i$  the eigenvalues of $\cM$,  then
\beqnas \tr U(X)U(Y)&=& \sum_i \frac{1}{(\lambda_i-X)(\lambda_i-Y)}= \frac{1}{Y-X}\sum_i \frac{1}{\lambda_i-Y}-\frac{1}{\lambda_i-X}\\&=&
\frac{1}{Y-X}\Big(\frac{P'(X)}{P(X)}-\frac{P'(Y)}{P(Y)}\Big).\eeqnas

For the formula for $\LL_{\bE,\Cl}(P)$, we start with
\beqnas \LL_{\bE,\Cl}(\log P)&= &\sum_{A,B,i,j} U(X)^{B,A}_{ji} \LL M^{A,B}_{ij} -\sum_{A,B,C,D,i,j,k,l}U(X)^{B,C}_{jk}U(X)^{D,A}_{l,i} \Gamma(\cM^{A,B}_{ij}, \cM^{C,D}_{kl})\\&=&
-\frac{1}{2}\sum_{A,B,C,D} E(A,B,C,D) U(X)^{B\DS C}_{jk}U(X)^{A\DS D}_{li}(\delta_{ik}\delta_{jl}+ (A\DS B|A\DS B)\delta_{il}\delta_{jk}),\eeqnas
where 
$$E(A,B,C,D)=\I_{A\DS B\DS C\DS D= \emptyset}(B\DS C|C)(A\DS D|A)(A\DS B|B)(C\DS D|D)= \I_{A\DS B\DS C\DS D= \emptyset}(A\DS C|A\DS C)$$
We obtain in the end
\beqnas \LL_{\bE,\Cl}(\log P)&=& -\frac{1}{2}\Bigg[\sum_{A,B,C}
(A\DS C|A\DS C)(B\DS C|B\DS C) \tr\big(U(X)^{B\DS C}\big)^2\\& &+ (B\DS C|A\DS C)(A\DS B|B\DS C)\big(\tr U(X)^{B\DS C}\big)^2\Bigg]\\
&=& -\frac{1}{2}\Big(\sum_{A} (A|A)\Big) \tr\big( U(X)^2\big)-2^{p-1} \sum_C (C|C) H(C)\big (\tr U(X)^C\big)^2.\eeqnas

\epf

If we are   interested in images of the Gaussian measure, we consider  the Ornstein-Uhlenbeck  operator  \[\LL_{\mathbb {OU},\Cl}(P)=\LL_{\bE,\Cl}(P)-D(P)\] where 
\[D=  \frac{1}{2}\Gamma(\|M \|^2,\cdot)\]
and
\[\|M \|^2=\sum_{i,j,A,B}(M^{A,B}_{i,j})^2\,.\]
If one is interested in images of the uniform measure on the unit sphere, we consider  instead the spherical operator
\beqna
\LL_{\bS,\Cl}\bS(P)=\LL_{\bE,\Cl}(P)-D^2(P)-(N-2)D(P), 
\eeqna
 where $N$ is the dimension on the Euclidean space in which the sphere is embedded, that is
 \[N=\frac{n2^p(n2^p+1)}{2}\,.\]
Observing their action on the characteristic polynomial, 
we have
\begin{eqnarray*}
D(\log P)=\frac{1}{2}\sum_{i,j,k,l,A,B,C,D}2M_{i,j}^{A,B}\partial_{M_{k,l}^{C,D}}(\log P)\Gamma_{\bE}(\cM^{A,B}_{ij}, \cM^{C,D}_{kl})
\\
=\frac{1}{2}\sum_{i,j,A\DS B\DS C\DS D= \emptyset}M_{i,j}^{A,B}(A\DS B|B\DS D)(U(X)^{D,C}_{j,i}+(A\DS B|A\DS B)U(X)^{D,C}_{i,j})
\\
=\sum_{A\DS B\DS C\DS D= \emptyset}(A\DS B|A\DS B)\tr\big(M^{A\DS B}U(X)^{A\DS B}\big)=2^p\tr\big(MU(X)\big)=2^p\big(n2^p-X\frac{P'}{P}\big)
\end{eqnarray*}
Hence
\beqna
D(P)=2^p(n2^pP-XP')
\eeqna
The carré du champ operator is the same for the Ornstein-Uhlenbeck operator than for the Laplace operator, whereas for the sphere, the carré du champ operator acting on the characteristic polynomial $P$ becomes  
\beqnas
&&\Gamma_{\bS,\Cl}(P(X),P(Y)) =\\
&&\frac{2^p}{Y-X}\big(P'(X)P(Y)-P'(Y)P(X)\big)-2^{2p}\big(n2^pP(X)-XP'(X)\big)\big(n2^pP(Y)-XP'(Y)\big)
\eeqnas

From Proposition~\ref{prop.gal.cliff}, one sees that the final expression depends on some specific factors for $\Cl(E)$ : the value of $\sum_A (A|A)$ and, for various $C\subset E$,  the value  of $H(C)= \sum_A (A|C)(C|A)$. We shall therefore restrict our attention to standard Clifford algebras, for which those computations may be explicitly done through some basic combinatorial arguments.

\section{ Symmetric matrices on standard Clifford algebras\label{sec.stand.cliff}}

Recall  that for standard Clifford algebras,  and with $E= \{1, \cdots, n\}$, one has, for any pair $(i,j)\in E^2$, $(i|j)= \I_{i<j}-\I_{j\leq i}$.

From $(A|B)= \prod_{i\in A, j\in B} (i|j)$, this immediately leads to
\beq\label{form.base.std.cliff}(A|A)= (-1)^{|A|(|A|+1)/2}, ~(A|B)(B|A)= (-1)^{|A||B|+ |A\cap B|}.\eeq

Notice also that for any $i\in E= \{1, \cdots, p\}$, $(i|E)= (-1)^i$.

\bprop\label{prop.combinat3} In a standard Clifford algebra with $|E|=p$, one has $$\sum_A (A|A)= \begin{cases} 2^{2m}(-1)^m & \hbox{ if $p= 4m$}\\
0&\hbox{ if $p= 4m+1$}\\
2^{2m+1}(-1)^{m+1} & \hbox{ if $p= 4m+2$}\\
2^{2m+2}(-1)^{m+1} &\hbox{ if $p= 4m+3$}\end{cases}$$

 \eprop
 \bpf From~\eqref{form.base.std.cliff}, one has
$$\sum_A (A|A)=\sum_k {p\choose k} (-1)^{k(k+1)/2}= \sum_k {p \choose 2k} (-1)^k - \sum_{k} {p \choose 2k+1} (-1)^{k}.$$
Comparing with
$$(1+i)^p=\sum_k {2k \choose p} (-1)^k + i\sum_{k} {2k+1 \choose p} (-1)^{k},$$ we see that 
$\sum_k {2k \choose p} (-1)^k $ is the real part of  $(1+i)^p$, while $\sum_{k} {2k+1 \choose p} (-1)^{k}$ is its imaginary part.
But $1+i= \sqrt{2}e^{i\pi/4}$, and therefore, for $p= 4m$,  
$(1+i)^{p}= 2^{2m}(-1)^m$,  for $p= 4m+1$, $(1+i)^{p}= (1+i)2^{2m}(-1)^m$, for  $p= 4m+2$, $(1+i)^p= i2^{2m+1}(-1)^m $ and for $p= 4m+3$, $(1+i)^{p}= (1+i)2^{2m}(-1)^m$. 
 It remains to collect the  various cases.
\epf

The following will also be useful

\bprop \label{prop.combinat} For a standard Clifford algebra $\Cl(E)$, and for  $B,C\subset E$, let 
\beqnas\begin{cases}S^e(B,C)= \sum_{A\subset C, |A|= 2k}  (A|B)(B|A),& S^e(B, \emptyset)= 1\\
 S^o(B,C)= \sum_{A\subset C, |A|= 2k+1} (A|B)(B|A), &S^o(B, \emptyset)= 0\end{cases}\eeqnas

Then, 
\beqnas\begin{cases}S^e(B,C)&= S^e( B\cap C,C),\\
 S^o(B,C)&= (-1)^{|B\cap C^c|}S^o(B\cap C,  C).\end{cases}\eeqnas
and, for $B\subset C$
 $$\begin{cases} S^e(B,C)= S^o(B,C)=0, & B\neq \emptyset, B\neq C\\
 S^e(\emptyset, C)= S^o(\emptyset, C)= 2^{|C|-1}, & C\neq \emptyset\\
  S^e(\emptyset, \emptyset)= 1, S^o(\emptyset, \emptyset)=0, & \\ 
 S^e(B,B)= -S^o(B,B)= 2^{|B|-1}&  |B|= 2k, B\neq \emptyset\\
 S^e(B,B)= S^o(B,B)= 2^{|B|-1}& |B|= 2k+1\end{cases}$$
 
 \eprop

\bpf For the first point,we decompose $B= (B\cap C) \cup (B\cap C^c)= B_1\cup B_2$. Then, if   $B_1\cap B_2 = \emptyset$
$$(A|B_1\cup B_2)(B_1\cup B_2|A)= (A|B_1)(B_1|A).(A|B_2)(B_2|A),$$  
and, if  $A\cap B_2= \emptyset$, then  $(A|B_2)(B_2|A)= (-1)^{|A||B_2|}$.

  It remains to study  $S^e(B,C)$ and $S^o(B,C)$ for $B\subset C$. Replacing $E$ by $C$, we are therefore bound to study the same quantity for  a standard Clifford algebra $\Cl(C)$.

Let us then fix  $C$ et $B\subset C$.  For  $C\neq \emptyset$
$S^e(\emptyset,C)= S^o(\emptyset, C)= 2^{|C|-1}$. When $B\neq \emptyset$ chose some point $i\in B$, and cut $\cP(C)$ into $\{A\subset C, x\in A\}$ and $\{A\subset C\setminus \{x\}\}$. 

Summing on $\cP(A)$, we get
\beqnas \begin{cases} S^e(B,C) &=  S^e(B\setminus x, C\setminus x) + (-1)^{|B|}S^o(B\setminus x, C\setminus x),\\ S^o(B)&=  (-1)^{|B|-1}S^e(B\setminus x, C\setminus x)-S^o(B\setminus x, C\setminus x).\end{cases}\eeqnas

In another way, setting  avec $U(B,C)= \bpm S^e(B,C)\\S^o(B,C)\epm$,  $U_k$ si $|B|= k$, on a 
$$U(B,C)= M_{\epsilon_k} U\big(B\setminus \{x\}, C\setminus \{x\}\big),$$  where $\epsilon_k= (-1)^{|B|}$ and 
\beq\label{def.matrices} M_1= \bpm1&1\\-1&-1\epm M_{-1}= \bpm 1&-1\\1&-1\epm.\eeq Setting  $S= \bpm 0&1\\1&0\epm$, which satisfies $S^2=1$, one has 
$$M_1M_{-1}=2(1-S), M_{-1}M_1= 2(1+S), M_1M_{-1}M_1 = 2M_1, M_{-1}M_1M_{-1}= 2M_{-1}, (1+S)M_{-1}= 2M_{-1}$$ from which  
$$(M_1M_{-1})^k= 2^{2k-1}(1-S), (M_{-1}M_1)^k= 2^{2k-1}(1+S),$$
and also
$$M_{-1}^2=0, M_1^2=0, M_1M_{-1}M_1=2M_1, M_{-1}M_1M_{-1}=2M_{-1}.$$ 
In the end, we get
$$\begin{cases}  |B|= 2k,&
U(B,C)= 2^{2k-1}(1+S) U(\emptyset, C\setminus B)\\ |B|= 2k+1,& U(B,C)= 2^{2k+1}(1-S)M_{-1}U(\emptyset, C\setminus B)\end{cases}.$$

Il remains to collect all the possible cases.

\epf

From  Proposition~\ref{prop.combinat}, one sees that, in a standard Clifford algebra $\Cl(E)$ with $|E|=n$, one has, for any $B\subset E$, with $B\neq \emptyset, E$, 
$$\sum_A (A|B)(B|A)= S^e(B,E)+S^o(B,E)=0.$$
Moreover,   $S^e(E,E)+S^o(E,E)= 0$ when  $|E|= 2k$. 
Therefore, $H(C)= 0$ unless $C= \emptyset$ or $C=E$, and $H(E)=0$ when $|E|= 2k$.
This leads to
\bprop \label{prop1}When $p= |E|$, 
if $\cM$ is $\Cl(E)$-symmetric, with  $P(X)=\det(M-XI)$ and $U(X)= (M-X\Id)^{-1}$, 
 $$\Gamma_{\bE,\Cl}\big(P(X),P(Y)\big)= \frac{2^p}{Y-X}\big(P'(X)P(Y)-P'(Y)P(X)\big),$$ and 

$$\begin{cases} \frac{\LL_{\bE,\Cl}(P)}{P} = \Gamma_{\bE,\Cl}(\log P)-2^{2m}(-1)^{m+1}(\tr U(X)^2) -\frac{1}{2} (\tr U(X))^2& \hbox{when  $|E|= 4m+2$}\\
\frac{\LL_{\bE,\Cl}(P)}{P} = \Gamma_{\bE,\Cl}(\log P)-2^{2m-1}(-1)^m(\tr U(X)^2) -\frac{1}{2} (\tr U(X))^2& \hbox{when $|E|= 4m$}\end{cases}$$
In particular, 
$$\frac{\LL_{\bE,\Cl}(P)}{P} = \begin{cases} (2^n+2^{2m}(-1)^m)\big( \frac{P'^2}{P^2}-\frac{P''}{P}\Big) -\frac{1}{2} \frac{P'^2}{P^2} & \hbox{ when  $|E|= 4m+2$}\\

(2^n+2^{2m-1}(-1)^{m+1})\big( \frac{P'^2}{P^2}-\frac{P''}{P}\Big) -\frac{1}{2} \frac{P'^2}{P^2} & \hbox{ when $|E|= 4m$}\end{cases} $$

As a consequence, one has $P(X)= Q(X)^a$, where $Q$ is a polynomial, where
$$\bcas a= 2^{4q},& \hbox{  when $n= 8q$}\\
 a= 2^{4q+2},& \hbox{  when $n= 8q+2$}\\
  a= 2^{4q+3},& \hbox{  when $n= 8q+4$}\\
   a= 2^{4q+3},& \hbox{  when $n= 8q+6$}\\
\ecas$$

Moreover, in those case,  for $\hat \LL= \frac{a}{2^{p}}\LL_{\bE,\Cl}$, $Q$ satisfies 
$$\hat \Gamma(Q(X),Q(Y))= \frac{1}{Y-X}\big(Q'(X)Q(Y)-Q'(Y)Q(X)\big)$$ and
$$\bcas \hat \LL(Q)= -\frac{1}{2} Q'', ~ p= 8q, ~p= 8q+6 & \hbox{~(real case)}\\
\hat \LL(Q)= -2 Q'', ~ p= 8q+2,~ p= 8q+4 & \hbox{~(quaternionic  case)}
\ecas
$$

\eprop

\bpf Using  Proposition~\ref{prop.gal.cliff}, the only term in the formula for $\LL_{\bE,\Cl} (P)$ which is not immediate to identify is $\tr U(X)^\emptyset$. But $\tr U(X)= 2^p\tr U(X)^\emptyset$, since only $U(X)^\emptyset$ appear in the diagonal blocs  of $U(X)$. Then, everything boils down to the computation of $\sum_A (A|A)$  given in Proposition~\ref{prop.combinat3}.
Then, the identification of $a$ comes from Equation~\eqref{eq.puiss} in Lemma~\ref{lem.puiss.P}. It is worth to observe that those equations have indeed integer roots in every  case. Moreover, Proposition~\ref{lem.puiss.P2} allows to assert that effectively, $P= Q^a$, and the rest is given again in Lemma~\ref{lem.puiss.P}.

\epf

It remains to deal with the case where $n$ is odd.  Then, the term $(\tr U(X)^E)^2$ appears in the formula for $\LL_{\bE,\Cl} (P)$.
But, when $(E|E)=-1$, then $U(X)^E$ is antisymmetric and $\tr U(X)^E=0$. Since 
$(E|E)= (-1)^{|E|(|E|+1)} $, this happens as soon as $|E|= 4m+1$.

This leads to

\bprop \label{prop2}Suppose that $p=|E|= 4m+1$.
Then, 
$$\Gamma_{\bE,\Cl}(P(X),P(Y))= \frac{2^p}{Y-X} \big(P'(X)P(Y)-P'(Y)P(X)\big),$$

$$\frac{\LL_{\bE,\Cl}(P)}{P}= 2^p(\frac{P'^2}{P^2}-\frac{P''}{P})-\frac{1}{2} \frac{P'^2}{P^2}.$$ 
Setting $a= 2^{2m+1}$ and  $P= Q^a$, and   for $\hat \LL= a2^{-p}\LL$, one has $\hat L (Q)= -Q''$. Then, $Q$ is a polynomial and the model corresponds to the complex case.

\eprop

It remains to deal with the case $p= 4m+3$, which turns out to be more delicate. Indeed, in those cases, $\Cl(E)$ is no longer simple, and splits into the direct sum of two ideals. From the Propositions~\ref{prop.gal.cliff} and~\ref{prop.combinat} we see that the set $E$ plays a special role in the analysis of $\LL_{\bE,\Cl}(P(X))$. 

 We already saw that in this situation, $(E|E)=1$, and for any $A\subset E$,
$(A|E)(E|A)= (-1)^{|A|(|E|+1)}= 1$, so that 
$\omega_E$ commutes to every element in the algebra and  satisfies $\omega_E^2=1$. Then, one may decompose the algebra $\Cl(E)$ into  the sum of the two ideals 
$\Cl(E)_+=\{ x\in \Cl(E), ~\omega_E x=x\}$ and $\Cl(E)_-=\{x \in \Cl(E), ~\omega_E x= -x\}$.
Symmetric matrices will also split into the direct sum of two symmetric matrices, and therefore the characteristic polynomial will be the product of characteristic polynomials.

We are therefore bound to consider separately the action on $\Cl(E)_+$ and $\Cl(E)_-$. We concentrate on the first one.
First observe that in this situation, for any $A\subset E$, $(A|E)= (E|A)= (A\DS E|E)= (E| A\DS E)$. From this, it is easy to see that 
$$\Cl(E)_+=\{X= \sum_A \lambda_A(\omega_A+ (A|E)\omega_{A\DS E})\}.$$

The action of the matrix $\sum_A M^A\omega_A$ on $\Cl(E)_+$ is the same as $\sum_A (A|E)M^{A\DS E} \omega_{A\DS E}$, and therefore  one may concentrate on matrices $\sum_A M^A \omega_A$ such that $M^{A\DS E} = (A|E) M^A$. This condition is clearly compatible with $(M^A)^t= (A|A) M^A$.  We therefore chose 
$$\Gamma_{\bE,\Cl}M_{ij}^A, M_{kl}^B)= \frac{1}{2}(\I_{A\DS B= \emptyset}+(A|E) \I_{A\DS B= E})(\delta_{ik}\delta_{jl}+ (A|A)\delta_{il}\delta_{jk}).$$

We may start the computation again, but it is simpler to observe that, setting $\sigma_A= \frac{1}{2}\big(\omega_A+ (A|E)\omega_{A\DS E}\big)$, one gets 
$\sigma_A\sigma_B = (A|B) \sigma_{A\DS B}$, and therefore  the  family $\sigma_A$ generates a standard Clifford algebra with with size $|E|-1$. Then, we boil down to a standard Clifford algebra with size $4m+2$,  and we see that we obtain the quaternionic case when $p= 8q+3$ and the real one when $p= 8q+7$. 

If one wants to describe the law of $P(X)$ in the case $p= 4m+3$, then we write $P= P_1P_2$, where $P_1$ and $P_2$ behave independently  as the previous ones.

We have such described all the laws of the spectra for symmetric matrices on standard Clifford algebras.

We thus recover Bott periodicity : in the following table, we give the algebra structure of $Cl(p)$, together with the dimension $d$ of the irreducible spaces in the third column,   the multiplicity $\alpha$  of the roots of the characteristic polynomial in the fourth, computed from the generator. In the last column, we indicate the parameter $a$ for which the law of the simple roots  $(\lambda_1< \cdots < \lambda_d)$ have density $\prod|\lambda_i-\lambda_j|^a$ with respect to the Lebesgue measure $d\lambda_1\cdots d\lambda_d$

$$\begin{tabular}{|c|c|c|c|c|}\hline
$|E|$& structure& $d$ &$\alpha$& a\\ \hline
\hline
Cl(1)&$\bC$&2&2&2\\\hline
Cl(2)&$\bH$&4&4&4\\\hline
Cl(3)&$\bH\oplus \bH$&4&4&4\\\hline
Cl(4)&$\bH[2]$& 8&8&4\\\hline
Cl(5)&$\bC[4]$& 8&8&2\\\hline 
Cl(6)&$\bR[8]$&8&8&1\\\hline
Cl(7)&$\bR[8]\oplus \bR[8]$&8&8&1\\\hline
Cl(8)&$\bR[16]$&16&16&1\\\hline
\end{tabular}
$$
then we tensorize by $\bR[16]$ through Bott's periodicity : $Cl(p+8)= \bR[16]\otimes Cl(p)$. (Here, $K[n]$ denotes the irreducible algebra of square  $n\times n$  matrices with coefficients in  the field  $K$). We may then observe that the multiplicity of the roots corresponds as expected to the dimension of the irreducible spaces, and that the parameter $a$ corresponds to the structure algebra : when the irreducible components are $K[n]$, then $a= 1,2, 4$ corresponding to the case where $K= \bR, \bC$ or $\bH$.

\brmq
 Considering the O-U operator $L_{\mathbb{OU},\Cl}$ described in the previous section, one gets here
 \beq
 L_{\mathbb{OU}}(P)=-CP''+(C-\frac{1}{2})\frac{P'^2}{P}-n2^{2p}P+2^pXP'
 \eeq
 and analogously to Lemma \ref{lem.puiss.P}, we get
 \[a\frac{L_{\mathbb{OU},\Cl}(P_1)}{P_1}=\frac{L_{\mathbb{OU},\Cl}(P)}{P}+\frac{1-a}{a}\Gamma_{\bE,\Cl}(\log P,\log P)\]
 which leads to
 \beq
 L_{\mathbb{OU},\Cl}(P_1)=-(C+2^p\frac{(1-a)}{a})P_1''+((C+2^p\frac{(1-a)}{a}-\frac{1}{2}a)\frac{P_1'^2}{P_1}-n\frac{2^{2p}}{a}P_1+2^pXP_1'
 \eeq 
 where $P=P_1^a$ and the constant $C$ differs according to $n$ ( see Propositions \eqref{prop1},\eqref{prop2}).
 Then, choosing $a$ as before, we can boil  down to the following relation:
 \beq
 L_{\mathbb{OU},\Cl}(P_1)=-\frac{a}{2}P_1''-n\frac{2^{2p}}{a}P_1+2^pXP_1'
 \eeq 
 
\ermq

\bibliographystyle{amsplain}

\providecommand{\bysame}{\leavevmode\hbox to3em{\hrulefill}\thinspace}
\providecommand{\MR}{\relax\ifhmode\unskip\space\fi MR }
\providecommand{\MRhref}[2]{%
  \href{http://www.ams.org/mathscinet-getitem?mr=#1}{#2}
}
\providecommand{\href}[2]{#2}
\begin{thebibliography}{}

\end{thebibliography}


\begin{thebibliography}{10}

\bibitem{AndGuioZeit}
G.W. Anderson, A.~Guionnet, and O.~Zeitouni.
\newblock An introduction to random matrices.
\newblock In {\em Cambridge Studies in Advanced Mathematics}, volume 118.
  Cambridge University Press, Cambridge, 2010.

\bibitem{atiyahbottshapiro}
M.~F. Atiyah, R.~Bott, and A.~Shapiro.
\newblock Clifford modules.
\newblock {\em Topology}, 3:3--38, 1964.

\bibitem{atiyahsinger}
M.~F. Atiyah and I.~M. Singer.
\newblock Index theory for skew-adjoint fredholm operators.
\newblock {\em Inst. Hautes {{\'E}}tudes Sci. Publ. Math.}, 37:5--26, 1969.

\bibitem{bglbook}
D.~Bakry, I.~Gentil, and M.~Ledoux.
\newblock {\em Analysis and {G}eometry of {M}arkov {D}iffusion {O}perators},
  volume 348 of {\em Grund. Math. Wiss.}
\newblock Springer, Berlin, 2013.

\bibitem{YanDoumerc}
Y.~Doumerc.
\newblock {\em Matrices al{{\'e}}atoires, processus stochastiques et groupes de
  r{\'e}flexions}.
\newblock PhD thesis, Universit\'e Toulouse 3, 2005.

\bibitem{dyson}
F.J. Dyson.
\newblock A brownian-motion model for the eigenvalues of a random matrix.
\newblock {\em Journ. of Mathematical Phys.}, 3:1191--1198, 1962.

\bibitem{erdosandco3}
L.~Erd{\H o}s, A.~Knowles, H.-T. Yau, and J.~Yin.
\newblock The local semicircle law for a general class of random matrices.
\newblock {\em Electronic Journ. of Prob.}, 18, 2013.

\bibitem{erdosandco1}
L.~Erd{\H o}s, S.~P{\'e}ch{\'e}, J.A. Ram{\'\i}rez, B.~Schlein, and H.-T. Yau.
\newblock Bulk universality for wigner matrices.
\newblock {\em Comm. Pure Appl. Math.}, 63:895--925, 2010.

\bibitem{erdosandco2}
L.~Erd{\H o}s, J.~Ram{\'\i}rez, B.~Schlein, T.~Tao, V.~Vu, and H.-T. Yau.
\newblock Bulk universality for wigner hermitian matrices with subexponential
  decay.
\newblock {\em Math. Res. Lett.}, 17:667--674, 2010.

\bibitem{Forrester}
P.J. Forrester.
\newblock Log--gases and random matrices.
\newblock In {\em London Mathematical Society Monographs Series}, volume~34.
  Princeton University Press, Princeton, 2010.

\bibitem{RamJag}
R.~Jagannathan.
\newblock On generalized clifford algebras and their physical applications.
\newblock In {\em The legacy of Alladi Ramakrishnan in the mathematical
  sciences}, pages 465--489. Springer, New York, 2010.

\bibitem{karoubi}
M.~Karoubi.
\newblock Alg{\`e}bres de clifford et {$K$}-th{\'e}orie.
\newblock {\em Ann. Sci. {\'E}cole Norm. Sup. (4)}, 1:161--270, 1968.

\bibitem{lounesto}
P.~Lounesto.
\newblock Clifford algebras and spinors (second edition).
\newblock In {\em London Mathematical Society Lecture Note Series}, volume 286.
  Cambridge University Press, Cambridge, 2001.

\bibitem{marcpastur}
L.~A. Mar{\v c}enko, V. A.;~Pastur.
\newblock Distribution of eigenvalues in certain sets of random matrices.
\newblock {\em Matematicheski{\u \i} Sbornik (N.S.)}, 72 (114):507--536, 1967.

\bibitem{Meh}
M.L. Mehta.
\newblock Random matrices.
\newblock In {\em Pure and Applied Mathematics (Amsterdam)}, volume 142.
  Elsevier/Academic Press, Amsterdam, 2004.

\bibitem{Sos99}
A.~Soshnikov.
\newblock Universality at the edge of the spectrum in wigner random matrices.
\newblock {\em Comm. in Math. Phys.}, 207:697--733, 1999.

\bibitem{Wig58}
E.P. Wigner.
\newblock On the distribution of the roots of certain symmetric matrices.
\newblock {\em Ann. Math}, 67:325--327, 1958.

\end{thebibliography}

\end{document}